\documentclass{article}

\usepackage[margin=1in]{geometry} 
\usepackage{amsmath,amsthm,amssymb}
\usepackage{graphicx}
\usepackage{color}
\usepackage{verbatim} 
\usepackage{url}

\usepackage{multirow}

\usepackage[dvipsnames]{xcolor}
\usepackage{subcaption}
\usepackage[export]{adjustbox}

\DeclareMathOperator*{\argmin}{\text{argmin}}
\graphicspath{{images/}}

\begin{document}

\title{Prediction of Discretization of GMsFEM using Deep Learning} 
\author{
Min Wang\thanks{Department of Mathematics, Texas A\&M University, College Station, TX 77843, USA (\texttt{wangmin@math.tamu.edu})}
\and
Siu Wun Cheung\thanks{Department of Mathematics, Texas A\&M University, College Station, TX 77843, USA (\texttt{tonycsw2905@math.tamu.edu})}
\and
Eric T. Chung\thanks{Department of Mathematics, The Chinese University of Hong Kong, Shatin, New Territories, Hong Kong SAR, China (\texttt{tschung@math.cuhk.edu.hk})}
\and
Yalchin Efendiev\thanks{Department of Mathematics \& Institute for Scientific Computation (ISC), Texas A\&M University,
College Station, Texas, USA (\texttt{efendiev@math.tamu.edu})}
\and
Wing Tat Leung\thanks{Institute for Computational Engineering and Sciences, 
The University of Texas at Austin, Austin, Texas, USA (\texttt{wleungo@ices.utexas.edu})}
\and
Yating Wang\thanks{Department of Mathematics, Texas A\&M University, College Station, TX 77843, USA (\texttt{wytgloria@math.tamu.edu})}
}

\maketitle

\begin{abstract}
In this paper, we propose a deep-learning-based approach to 
a class of multiscale problems. 
THe Generalized Multiscale Finite Element Method (GMsFEM) 
has been proven successful as a model reduction technique of 
flow problems in heterogeneous and high-contrast porous media. 
The key ingredients of GMsFEM include mutlsicale basis functions and 
coarse-scale parameters, which are obtained from solving 
local problems in each coarse neighborhood. 
Given a fixed medium, these quantities are precomputed by solving 
local problems in an offline stage, 
and result in a reduced-order model. 
However, these quantities have to be re-computed in case of varying media. 
The objective of our work is to make use of deep learning techniques 
to mimic the nonlinear relation between the permeability field 
and the GMsFEM discretizations, and use neural networks to perform 
fast computation of GMsFEM ingredients repeatedly for a class of media. 
We provide numerical experiments to investigate the predictive power of neural networks 
and the usefulness of the resultant multiscale model in solving channelized porous media flow problems. 
\end{abstract}

\section{Introduction}
Multiscale features widely exist in many engineering problems. 
For instance, in porous media flow, the media properties typically 
vary over many scales and contain high contrast. 
Multiscale Finite Element Methods (MsFEM) \cite{eh09,hw97,jennylt03} 
and Generalized Multiscale Finite Element Methods (GMsFEM) \cite{GMsFEM13,review16} 
are designed for solving multiscale problems using local model reduction techniques. 
In these methods, the computational domain is partitioned into 
a coarse grid $\mathcal{T}^H$, which does not necessarily resolve
all multiscale features. We further perform a refinement of $\mathcal{T}^H$ 
to obtain a fine grid $\mathcal{T}^h$, which essentially resolves all multiscale features. 
The idea of local model reduction in these methods is based on 
idenfications of local multiscale basis functions supported in coarse regions 
on the fine grid, and replacement of the macroscopic equations 
by a coarse-scale system using a limited number of local multiscale basis functions. 
As in many model reduction techniques, the computations of multiscale basis functions, 
which constitute a small dimensional subspace, 
can be performed in an offline stage. 
For a fixed medium, these multiscale basis functions are 
reusable for any force terms and boundary conditions. 
Therefore, these methods provide a substantial computational savings 
in the online stage, in which a coarse-scale system is constructed and 
solved on the reduced-order space. 

However, difficulties arise in situations with uncertainties in the media properties in some local regions, 
which are common for oil reservoirs or aquifers. 
One straightforward approach for quantifying the uncertainties 
is to sample realizations of media properties. 
In such cases, it is challenging to find an offline principal component subspace 
which is able to universally solve the multiscale problems with different media properties. 
The computation of multiscale basis functions has to be performed in an online procedure 
for each medium. Even though the multiscale basis functions are 
reusable for different force terms and boundary conditions, 
the computational effort can grow very huge for a large number of realizations of media properties. 
To this end, building a functional relationship between the media properties and 
the multiscale model in an offline stage can avoid repeating expensive computations 
and thus vastly reduce the computational complexity. 
Due to the diversity of complexity of the media properties, 
the functional relationship is highly nonlinear. 
Modelling such a nonlinear functional relationship typically
involves high-order approximations. 
Therefore, it is natural to use machine learning techniques 
to devise such complex models. 
In \cite{bayesian17,bayesian18}, the authors make use of a Bayesian approach 
for learning multiscale models and incorporating essential observation data 
in the presence of uncertainties.  

Deep neural networks is one class of machine learning 
algorithm that is based on an artificial neural network, 
which is composed of a relatively large number of layers 
of nonlinear processing units, called neurons, 
for feature extraction. The neurons are connected to other neurons 
in the successive layers. The information propagates from the input, 
through the intermediate hidden layers, and to the output layer. 
In the propagation process, the output in each layer 
is used as input in the consecutive layer. 
Each layer transforms its input data into a 
little more abstract feature representation. 
In between layers, a nonlinear activation 
function is used as the nonlinear transformation on the input, 
which increases the expressive power of neural networks. 
Recently, deep neural network (DNN) has been successfully used to 
interpret complicated data sets and applied to tasks with pattern recognition, such as 
image recognition, speech recognition and natural language processing
\cite{dcnn_im2012, dnn_speech_2012, HK_resnet2016}. 
Extensive researches have also been conducted on investigating 
the expression power of deep neural networks 
 \cite{Cybenko1989, Hornik1991, Csaji2001, Telgrasky2016, Poggio2016, Hanin2017}. 
Results show that neural networks can represent and approximate
a large class of functions.
Recently, deep learning has been applied to model reductions and 
partial differential equations. 
In \cite{deepconv_Zabaras}, the authors studied 
deep convolution networks for surrogate model construction. 
on dynamic flow problems in heterogeneous media. 
In \cite{Shi_resnet}, the authors studied the relationship 
between residual networks (ResNet) and characteristic equations of linear transport, 
and proposed an interpretation of deep neural networks 
by continuous flow models. 
In \cite{E_deepRitz}, the authors combined the idea of 
the Ritz method and deep learning techniques to
solve elliptic problems and eigenvalue problems. 
In \cite{Ying_paraPDE}, a neural network has been designed  
to learn the physical quantities of interest as a function of 
random input coefficients. 
The concept of using deep learning to generate a reduced-order model for a dynamic flow 
has been applied to proper orthogonal decomposition (POD) global model reduction \cite{DNN_POD} 
and nonlocal multi-continuum upscaling (NLMC) \cite{DNN_NLMC}. 

In this work, we propose a deep-learning-based method 
for fast computation of the GMsFEM discretization. 
Our approach makes use of deep neural networks as a fast proxy to compute 
GMsFEM discretizations for flow problems in channelized porous media with uncertainties. 
More specifically, neural networks are used to express 
the functional relationship between the media properties and 
the multiscale model. 
Such networks are built up in an offline stage. 
Sufficient sample pairs are required to ensure the 
expressive power of the networks. 
With different realizations of media properties, one can 
use the built network and avoid computations of local problems and spectral problems. 

The paper is organized as follows. 
We start with the underlying partial differential equation 
that describes the flow within a heterogeneous media and 
the main ingredients of GMsFEM in Section~\ref{sec:prelim}. 
Next, in Section~\ref{sec:dl}, we present the idea of using deep learning 
as a proxy for prediction of GMsFEM discretizations. 
The networks will be precisely defined and 
the sampling will be explained in detail. 
In Section~\ref{sec:num}, we present numerical experiments to show the effectiveness of 
our presented networks on several examples with different configurations. 
Finally, a conclusive discussion is provided in Section~\ref{sec:conclusion}.

\section{Preliminaries}
\label{sec:prelim}

In this paper, we are considering the flow problem in highly heterogeneous media 
\begin{equation}\label{eq:CGonline_original}
	\begin{split}
	-\text{div} \big(\kappa \, \nabla u  \big)  =f \quad &\text{in} \quad \Omega,\\
	u  = 0 \quad \text{or} \quad  \frac{\partial u}{\partial n} = 0  \quad &\text{on} \quad \partial \Omega,
	\end{split}
\end{equation}
where $\Omega$ is the computational domain, 
$\kappa$ is the permeability coefficient in $L^{\infty}(\Omega)$, and 
$f$ is a source function in $L^2(\Omega)$.
We assume the coefficient $\kappa$ is highly heterogeneous with high contrast. 
The classical finite element method for solving \eqref{eq:CGonline_original} numerically is given by: 
find $u_h \in V_h$ such that
\begin{equation}\label{eq:FEM_variational}
\begin{split}
a(u_h,v) = \int_{\Omega}\kappa\nabla u_h\cdot\nabla v ~dx &= \int_{\Omega} fv ~dx = (f,v) \quad \text{for all} \, \, v\in V_{h},\\
\end{split}
\end{equation}
where $V_h$ is a standard conforming finite element space over a partition 
$\mathcal{T}_h$ of $\Omega$ with mesh size $h$. 

However, with the highly heterougeneous property of coefficient $\kappa$, the mesh size $h$ has to be taken extremely small to capture the underlying fine-scale features of $\kappa$. This ends up with a large computational cost. 
Generalized Multiscale Finite Element Method (GMsFEM) \cite{GMsFEM13,review16} serves as a model reduction technique to 
reduce the number of degree of freedom and 
attain both efficiency and accuracy to a considerable extent. 
GMsFEM has been successfully extended to other formulations and applied to other problems. 
Here we provide a brief introduction of the main ingredients of GMsFEM. 
For a more detailed discussion of GMsFEM and related concepts, the reader is referred to 
\cite{WaveGMsFEM,MixedGMsFEM,eglp13oversampling,randomized2014,chung2015residual}. 

In GMsFEM, we define a coarse mesh $\mathcal{T}^H$ over the domain $\Omega$ and 
refine to obtain a fine mesh $\mathcal{T}^h$ with mesh size $h \ll H$, 
which is fine enough to restore the multiscale properties of the problem. 
Multiscale basis functions are defined on coarse grid blocks using linear combinations of 
finite element basis functions on $\mathcal{T}^h$, and 
designed to resolve the local multiscale behaviors of the exact solution. 
The multiscale finite element space $V_{\text{ms}}$, 
which is a principal component subspace of the conforming finite space $V_h$ with 
$\text{dim}(V_{\text{ms}})\ll \text{dim}(V_{h})$, 
is constructed by the linear span of multiscale basis functions. 
The multiscale solution $u_{\text{ms}}\in V_{\text{ms}}$ is then defined by 
\begin{equation}
a(u_{\text{ms}},v) = (f,v) \quad \text{for all} \, \, v\in V_{\text{ms}}.
\end{equation}

\begin{figure}[ht]
	\centering
	\includegraphics[width = 0.8\textwidth]{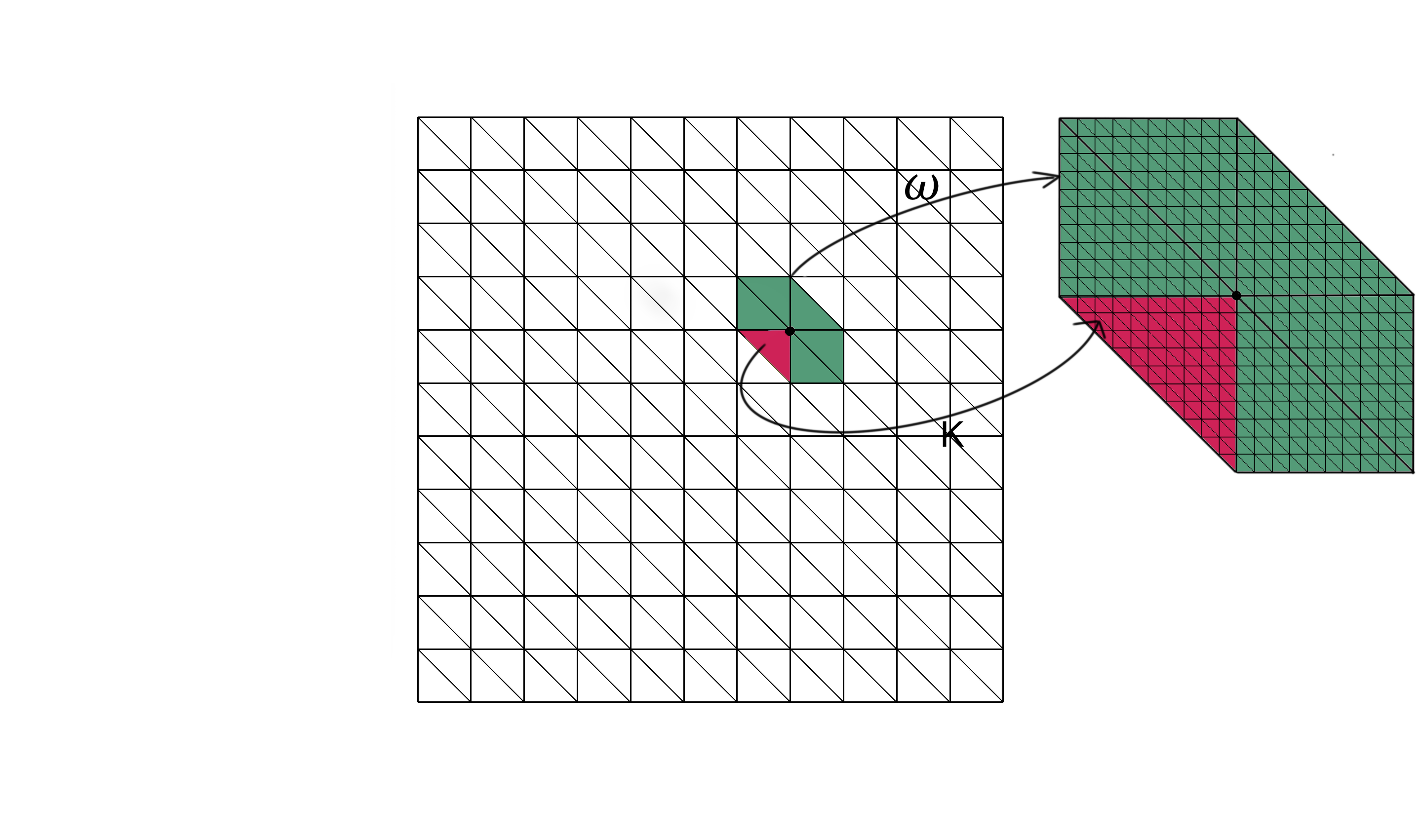}
	\caption{An illustration of coarse mesh (left), a coarse neighborhood and coarse blocks (right)}
	\label{fig:mesh}
\end{figure}

In this work, we consider the identification of dominant modes for solving \eqref{eq:CGonline_original} 
by multiscale basis functions, including spectral basis functions and simplified basis functions, in GMsFEM. 
Here, we present the details of the construction of multiscale basis functions in GMsFEM. 
Let $N_x=\{x_i\; | 1\leq i\leq N_v\}$ be the set of nodes of the coarse mesh $\mathcal{T}^H$. 
For each coarse grid node $x_i \in N_x$, the coarse neighborhood $\omega_i$ is defined by
\begin{equation}
\omega_i=\bigcup\{ K_j\in\mathcal{T}^H; ~~~ x_i\in \overline{K}_j\},
\end{equation}
that is, the union of the coarse elements $K_j \in\mathcal{T}^H$ containing the coarse grid node $x_i$. 
An example of the coarse and fine mesh, coarse blocks and a coarse neighborhood is shown in Figure~\ref{fig:mesh}. 
For each coarse neighbourhood $\omega_i$, 
we construct multiscale basis functions $\{\phi_j^{\omega_i}\}_{j=1}^{L_i}$ 
supported on $\omega_i$. 

For the construction of spectral basis functions, 
we first construct a snapshot space $V_{\text{snap}}^{(i)}$ 
spanned by local snapshot basis functions $\psi_{\text{snap}}^{i,k}$ 
for each local coarse neighborhood $\omega_i$. 
The snapshot basis function $\psi_{\text{snap}}^{i,k}$ is the solution of a local problem 
\begin{equation}\label{eq:harmonic}
\begin{split}
	-\text{div} (\kappa \nabla \psi_{\text{snap}}^{i,k}) = 0, &\quad \text{in} \quad \omega_i \\
	\psi_{\text{snap}}^{i,k} = \delta_{i,k}, & \quad \text{on} \quad \partial \omega_i
\end{split}	
\end{equation}
The fine grid functions $\delta_{i,k}$ is a function defined for all $x_s \in \partial \omega_i$, 
where $\{x_s\}$ denote the fine degrees of freedom on the boundary of local coarse region $\omega_i$. In specific, 
\[
\delta_{i,k}(x_s) =
\begin{cases}
1 & \text{if $s = k$} \\
0 & \text{if $s \neq k$} 
\end{cases}
\]
The linear span of these harmonic extensions forms the local snapshot space 
$V_{\text{snap}}^{(i)} = \text{span}_{k} \{\psi_{\text{snap}}^{i,k}\}$. 
One can also use randomized boundary conditions to reduce the computational cost 
associated with snapshot calculations \cite{randomized2014}. 
Next, a spectral problem is designed based on our analysis and 
used to reduce the dimension of local multiscale space. 
More precisely, we seek for eigenvalues $\lambda_m^{i}$ and 
corresponding eigenfunctions $\phi^{\omega_i}_m \in V_{\text{snap}}^{(i)}$ satisfying
\begin{equation}\label{eq:spectral}
a_i(\phi^{\omega_i}_m, v) = \lambda_m^{i} s_i(\phi^{\omega_i}_m, v), \quad \forall v \in V_{\text{snap}}^{(i)}, 
\end{equation}
where the bilinear forms in the spectral problem are defined as
\begin{equation}
\begin{split}
a_i(u, v) & = \int_{\omega_i} \kappa \nabla u \cdot \nabla v, \\
s_i(u, v) & = \int_{\omega_i}\tilde{\kappa} u  v, 
\end{split}
\end{equation}
where $\tilde{\kappa}  = \sum_{j} \kappa |\nabla \chi_j|^2$, 
and $\chi_j$ denotes the multiscale partition of unity function. 
We arrange the eigenvalues $\lambda_m^i$ of the spectral problem \eqref{eq:spectral} in ascending order, 
and select the first $l_i$ eigenfunctions $\{ \phi^{\omega_i}_m \}_{m=1}^{l_i}$ 
corresponding to the small eigenvalues 
as the multiscale basis functions.

An alternative way to construct the multiscale basis function 
is using the idea of simplified basis functions. 
This approach assumes the number of channels and 
position of the channalized permeability field are known. 
Therefore we can obtain multiscale basis functions $\{\phi_m^{\omega_i}\}_{m=1}^{l_i}$
using these information and without solving the spectral problem \cite{NLMC}.

Once the multiscale basis functions are constructed, 
the span of the multiscale basis functions will form the offline space
\begin{equation}
\begin{split}
V_{ms}^{(i)} & = \text{span}\{\phi_m^{\omega_i}\}_{m=1}^{l_i}, \\
V_{ms} & = \oplus_i V_{ms}^{(i)}.
\end{split}
\end{equation}
We will then seek a multlscale solution $u_{ms} \in V_{ms}$ satisfying 
\begin{equation}
a(u_{ms}, v) = (f,v) \quad \text{ for all } v \in V_{ms}.
\label{eq:var_ms}
\end{equation}
which is a Galerkin projection of the \eqref{eq:CGonline_original} onto $V_{ms}$, 
and can be written as a system of linear equations 
\begin{equation}
A_c u_c = b_c,
\label{eq:coarse_system}
\end{equation} 
where $A_c$ and $b_c$ are the coarse-scale stiffness matrix and load vector. 
If we collect all the multiscale basis functions and 
arrange the fine-scale coordinate representation 
in columns, we obtain the downscaling operator $R$. 
Then the fine-scale representation of the multiscale solution is given by
\begin{equation}
u_{ms} = R u_c.
\label{eq:downscale}
\end{equation}

\section{Deep Learning for GMsFEM}\label{sec:dl}

\begin{figure}[ht]
  \centering
  \includegraphics[width=.8\linewidth]{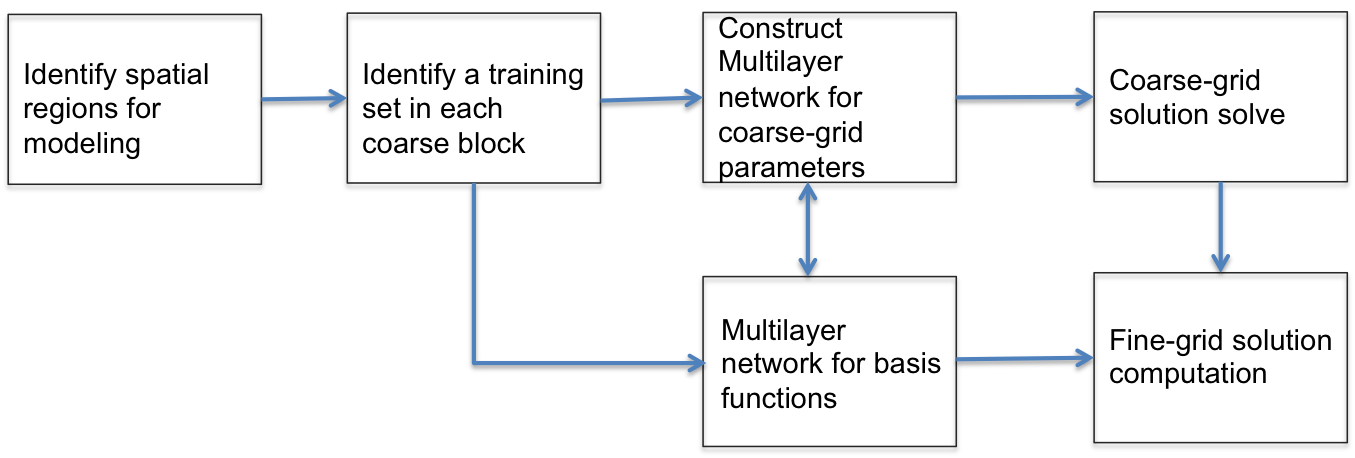}
  \caption{A flow chart in illustrating the idea of using deep learning in the GMsFEM framework. }
  \label{fig:chart1}
\end{figure}

In real world applications, there are uncertainties within some local regions 
of the permeability field $\kappa$ in the flow problem. 
Thousands of forward simulations are needed to quantify the uncertainties of the flow solution. 
GMsFEM provides us with a fast solver to compute the solutions accurately and efficiently. 
Considering that there is a large amount of simulation data, 
we are interested in developing a method utilizing the existing offline data 
and reducing direct computational effort later. 
In this work, we aim at using DNN to model the relationship between heterogeneous 
permeability coefficient $\kappa$ and the key ingredients of GMsFEM solver, 
i.e., coarse scale stiffness matrices and multiscale basis functions. 
When the relation is built up, we can feed the network any 
realization of the permeability field and obtain the corresponding GMsFEM ingredients, 
and further restore fine-grid GMsFEM solution of \eqref{eq:CGonline_original}.
The general idea of utilizing deep learning in the GMsFEM framework 
is illustrated in Figure~\ref{fig:chart1}.
%
%
%
%
%

\begin{figure}[ht]
  \centering
  \includegraphics[width=.8\linewidth]{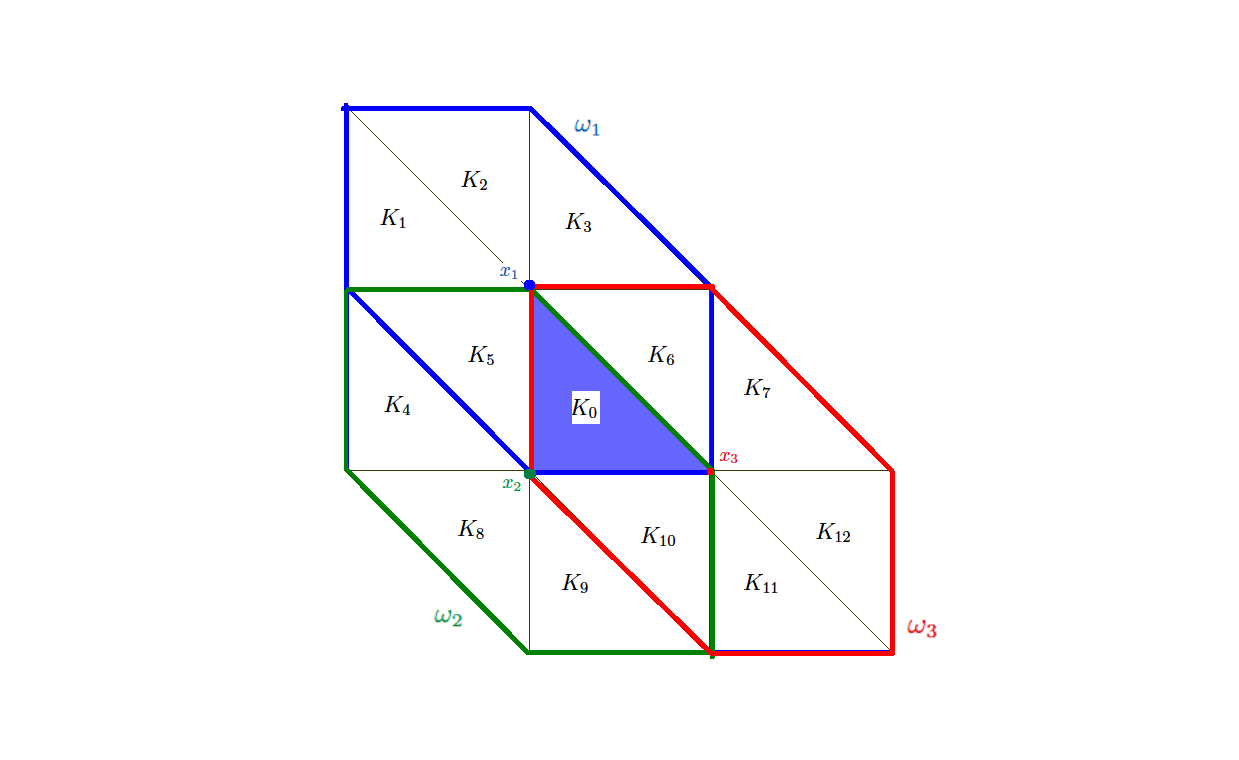}
  \caption{An illustration of a target coarse block $K_0$ and related neighborhoods.}
  \label{fig:kandomega}
\end{figure}

Suppose that there are uncertainties for the heterogeneous coefficient 
in a local coarse block $K_0$, which we call the target block, 
and the permeability outside the target block remains the same. 
For example, for a channelized permeability field, 
the position, location and the permeability values of the channels in the target block can vary. 
The target block $K_0$ is inside $3$ coarse neighborhoods, 
denoted by $\omega_1, \omega_2, \omega_3$. 
The union of the $3$ neighborhoods, i.e. 
\[ \omega^+(K_0) = \omega_1 \cup \omega_2 \cup \omega_3, \]
are constituted of by the target block $K_0$ and $12$ other coarse blocks, 
denoted by $\{K_l\}_{l=1}^{12}$
A target block and its surrounding neighborhoods 
are depicted in Figure~\ref{fig:kandomega}. 

For a fixed permeability field $\kappa$, 
one can compute the multiscale basis functions $\phi_m^{\omega_i}(\kappa)$ 
defined by \eqref{eq:spectral}, for $i = 1, 2, 3$,
and the local coarse-scale stiffness matrices $A_c^{K_l}(\kappa)$, defined by
\begin{equation}\label{eq:mat_element}
[A_c^{K_l}(\kappa)]^{i,j}_{m,n} = \int_{K_l} \kappa \nabla \phi_m^{\omega_i}(\kappa) \cdot \nabla \phi_n^{\omega_j}(\kappa).
\end{equation}
for $l = 0, 1, \ldots, 12$.
We are interested in constructing the maps $g_B^{m,i}$ and $g_M^l$, where
\begin{itemize}
\item $g_B^{m,i}$ maps the permeability coefficient $\kappa$ 
to a local multiscale basis function $\phi^{\omega_i}_m$, 
where $i$ denotes the index of the coarse block, and 
$m$ denotes the index of the basis in coarse block $\omega_i$
\[
g_B^{m,i}: \kappa \mapsto \phi_m^{\omega_i}(\kappa),
\]
\item $g_M^l$ maps the permeability coefficient $\kappa$ to 
the coarse grid parameters $A_c^{K_l}$  ($l = 0, \cdots, 12$)
 \[
g_M^l: \kappa \mapsto A_c^{K_l}(\kappa).
\]
\end{itemize}
In this work, our goal is to make use of deep learning to 
build fast approximations of these quantities 
associated with the uncertainties in the permeability field $\kappa$, 
which can provide fast and accurate solutions to the 
heterogeneous flow problem \eqref{eq:CGonline_original}.

For each realization $\kappa$, 
one can compute the images of $\kappa$ 
under the local multiscale basis maps $g_B^{m,i}$ 
and the local coarse-scale matrix maps $g_M^l$. 
These forward calculations serve as training samples for 
building a deep neural network for approximation 
of the corresponding maps, i.e. 
\begin{equation}
\begin{split}
\mathcal{N}_B^{m,i}(\kappa) & \approx g_B^{m,i}(\kappa), \\
\mathcal{N}_M^l(\kappa) & \approx g_M^l(\kappa).
\end{split}
\end{equation}
In our networks, the permeability field $\kappa$ is the input, 
while the multiscale basis functions $\phi^{\omega_i}_m$ 
and the coarse-scale matrix $A_c^{K_l}$ are the outputs.
Once the neural networks are built, 
we can use the networks to compute the 
multiscale basis functions and coarse-scale parameters 
in the associated region for any permeability field $\kappa$. 
Using these local information from the neural networks  
together with the global information which can be pre-computed, 
we can form the downscale operator $R$ with the 
multiscale basis functions, 
form and solve the linear system \eqref{eq:coarse_system}, 
and obtain the multiscale solution by \eqref{eq:downscale}. 

\subsection{Network architecture}

\begin{figure}[ht]
  \centering
  \includegraphics[width=.6\linewidth]{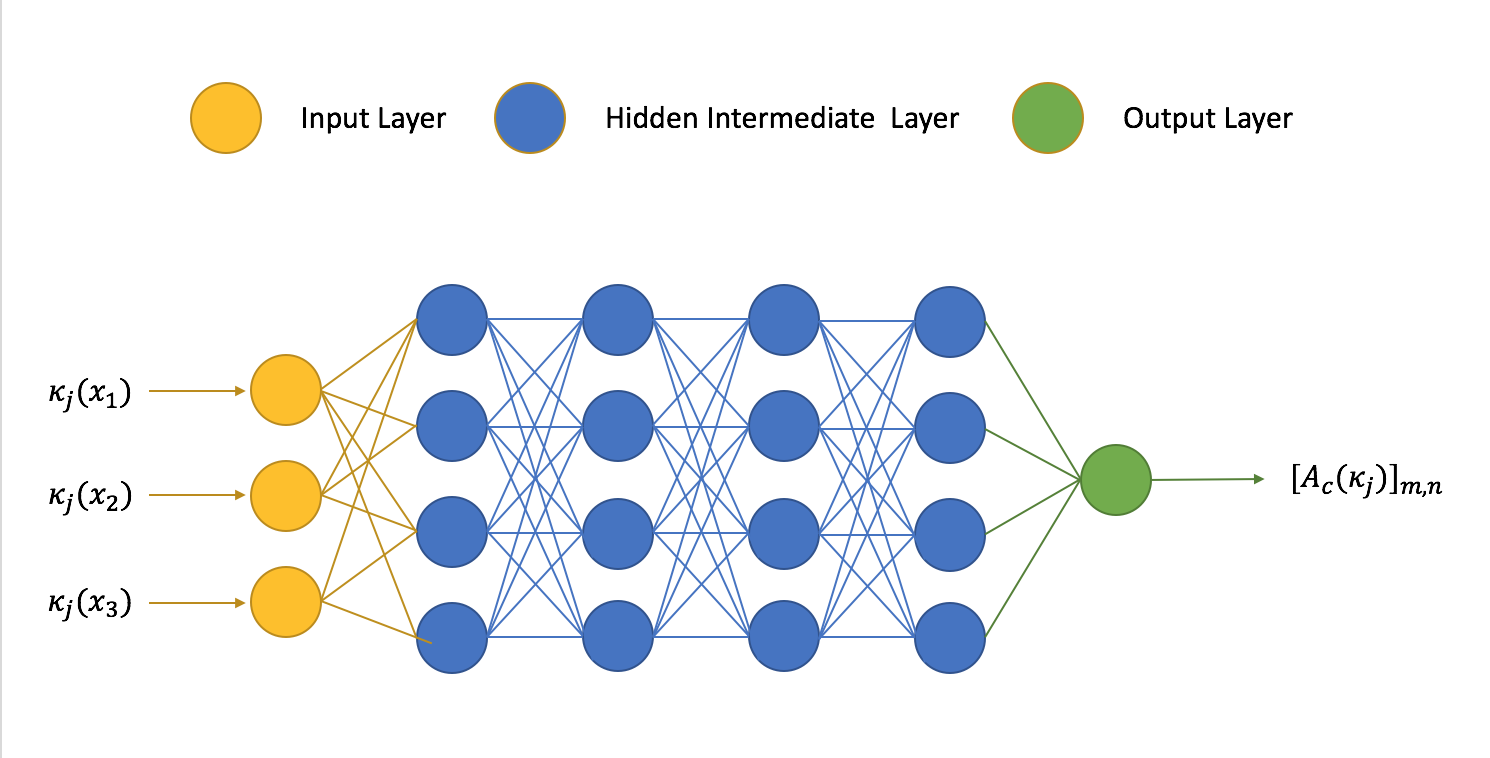}
  \caption{An illustration of a deep neural network.}
  \label{fig:dnn_illustration}
\end{figure}

In general, a $L$-layer neural network $\mathcal{N}$ can be written in the form 
\begin{equation*}
\mathcal{N}(x; \theta) = \sigma(W_L \sigma (\cdots  
\sigma(W_2  \sigma(W_1 x + b_1) + b_2) \cdots   ) + b_L),
\end{equation*} 
where $\theta : = (W_1, W_2, \cdots, W_L, b_1, b_2, \cdots, b_L)$, 
$W$'s are the weight matrices and $b$'s are the bias vectors, 
$\sigma$ is the activation function, 
$x$ is the input. 
Such a network is used to approximate the output $y$.
Our goal is then to find $\theta^*$ by solving an optimization problem
\begin{equation*}
\theta^{*} = \argmin_{\theta} ~~\mathcal{L}(\theta),
\end{equation*}
where $\mathcal{L}(\theta)$ is called loss function, 
which measures the mismatch between the image of the input $x$ 
under the the neural network $\mathcal{N}(x,y;\theta)$ and the desired output $y$ 
in a set of training samples $(x_j, y_j)$.
In this paper, we use the mean-squared error metric to be our loss function
\begin{equation*}
\mathcal{L}(\theta) =  \frac{1}{N}\sum_{j=1}^{N} \|y_j - \mathcal{N}(x_j; \theta) \|^2_2,
\end{equation*}
where $N$ is the number of the training samples. 
An illustration of a deep neural network is 
shown in Figure~\ref{fig:dnn_illustration}.

Suppose we have a set of different realizations of the permeability 
$\{\kappa_1, \kappa_2, \cdots, \kappa_N\}$ in the target block. 
In our network, the input $x_j = \kappa_j \in \mathbb{R}^d$ 
is a vector containing the permeability image pixels in the target block. 
The output $y_j$ is an entry of the local stiffness matrix $A_c^{K_l}$, 
or the coordinate representation of a multiscale basis function $\phi_m^{\omega_i}$. 
We will make use of these sample pairs $(x_j, y_j)$ to train a deep neural network 
$\mathcal{N}_B^{m,i}(x; \theta_B^*)$ and $\mathcal{N}_M^l(x; \theta_M^*)$ 
by minimizing the loss function with respect to the network parameter $\theta$, 
such that the trained neural networks can approximate the functions $g_B^{m,i}$ and $g_M^l$, respectively.
Once the neural is constructed, for some given new permeability coefficient $\kappa_{{N+1}}$, 
we use our trained networks to compute a fast prediction of the outputs, 
i.e. local multiscale basis functions $\phi_m^{\omega_i,\text{pred}}$ by 
\[ \phi_m^{\omega_i,\text{pred}}(\kappa_{N+1}) = 
\mathcal{N}_B^{m,i}(\kappa_{{N+1}};\theta_B^*) \approx 
g_B^{m,i}(\kappa_{N+1}) = \phi_m^{\omega_i}(\kappa_{N+1}), \]
and local coarse-scale stiffness matrix $A_c^{K_l,\text{pred}}$ by 
\[ A_c^{K_l,\text{pred}}(\kappa_{N+1}) =
 \mathcal{N}_M^l(\kappa_{{N+1}};\theta_M^*)
\approx g_M^l(\kappa_{N+1}) = A_c^{K_l}(\kappa_{N+1}). \]





\subsection{Network-based multiscale solver}
Once the neural networks are built, we can 
assemble the predicted multiscale basis functions 
to obtain a prediction 
$R^{\text{pred}}$ for the downscaling operator, 
and assemble the predicted local coarse-scale stiffness matrix $A_c^{K_l,\text{pred}}$
in the global matrix $A_c^{\text{pred}}$. 
Following \eqref{eq:coarse_system} and \eqref{eq:downscale},
we solve the predicted coarse-scale coefficient vector 
$u_c^{\text{pred}}$ from the following linear system
\begin{equation}\label{eq:coarse_system_pred}
A_c^{\text{pred}} u_c^{\text{pred}} = b_c,
\end{equation}
and obtain the predicted multiscale solution $u_{ms}^{\text{pred}}$ by
\begin{equation}
u_{ms}^{\text{pred}} = R^{\text{pred}} u_c^{\text{pred}}.
\label{eq:downscale_pred}
\end{equation}

\section{Numerical Results}
\label{sec:num}

%
%
%

In this section, we present some numerical results for 
predicting the GMsFEM ingredients and solutions 
using our proposed method. 
We consider permeability fields $\kappa$ with high-contrast channels 
inside the domain $\Omega = (0,1)^2$, 
which consist of uncertainties in a target cell $K_0$. 
More precisely, we consider a number of random realizations of permeability fields 
$\kappa_1,\kappa_2, \kappa_3,\cdots, \kappa_{N+M}$. 
Each permeability field contains two high-conductivity channels, 
and the fields differ in the target cell $K_0$ by:
\begin{itemize}
\item in Experiment 1, the channel configurations are all distinct, 
and the permeability coefficients inside the channels are fixed in each sample 
(see Figure~\ref{fig:kappa_one_samples} for illustrations), and
\item in Experiment 2, the channel configurations are randomly chosen among 5 configurations, and 
the permeability coefficients inside the channels follow a random distribution
(see Figure~\ref{fig:sine_samples} for illustrations).
\end{itemize}
In these numerical experiments, we assume there are uncertainties in only the 
target block $K_0$. The permeability field in $\Omega \setminus K_0$ 
is fixed across all the samples. 

We follow the procedures in Section~\ref{sec:dl} and
to generate sample pairs using GMsFEM. 
Local spectral problems are solved to obtain the multiscale basis functions $\phi_m^{\omega_i}$. 
In the neural network, the permeability field $x = \kappa$ is considered to be the input,
while the local multiscale basis functions $y = \phi_m^{\omega_i}$ and 
local coarse-scale matrices $y = A_c^{K_l}$ are regarded as the output. 
These sample pairs are divided into the training set and the learning set in a random manner. 
A large number $N$ of realizations, namely $\kappa_1, \kappa_2, \ldots, \kappa_N$, 
are used to generate sample pairs in the training set, 
while the remaining $M$ realizations, namely, $\kappa_{N+1}, \kappa_{N+2}, \ldots, \kappa_{N+M}$ 
are used in testing the predictive power of the trained network. 
We remark that, for each basis function and each local matrix, 
we solve an optimization problem in 
minimizing the loss function defined by the sample pairs in the training set, 
and build a separate deep neural network. 
We summarize the network architectures for 
training local coarse scale stiffness matrix and multiscale basis functions as below:
\begin{itemize}
\item For the multiscale basis function $\phi_m^{\omega_i}$, 
we build a network $\mathcal{N}_B^{m,i}$ using
\begin{itemize}
\item {Input}: Vectorized permeability pixels values $\kappa$,
\item {Output}: Coefficient vector of multiscale basis $\phi_m^{\omega_i}(\kappa)$ on coarse neighborhood $\omega_i$,
\item {Loss Function}: Mean squared error $\displaystyle\frac{1}{N}\sum_{j=1}^{N} ||\phi_m^{\omega_i}(\kappa_j)- \mathcal{N}_B^{m,i} (\kappa_j;\theta_B)||_2 ^2$,
\item {Activation Function}: Leaky ReLu function,
\item {DNN structure}: 10-20 hidden layers, each layer have 250-350 neurons,
\item {Training Optimizer}: Adamax.
\end{itemize}
\item For the local coarse scale stiffness matrix $A_c^{K_l}$, 
we build a network $\mathcal{N}_M^l$ using
\begin{itemize}
\item {Input}: Vectorized permeability pixels values $\kappa$, 
\item {Output}: Vectorized coarse scale stiffness matrix $A_c^{K_l}(\kappa)$ on the coarse block $K_l$,
\item {Loss Function}: Mean squared error $\displaystyle\frac{1}{N}\sum_{j=1}^{N} ||A_c^{K_l}(\kappa_j) - \mathcal{N}_M^l (\kappa_j;\theta_M)||_2 ^2$,
\item {Activation Function}: ReLu function (Rectifier),
\item {DNN structure}: 10-16 hidden layers, each layer have 100-500 neurons,
\item {Training Optimizer}: Proximal Adagrad.
\end{itemize}
\end{itemize}
For simplicity, the activation functions ReLU function \cite{glorot11} and Leaky ReLU function are used 
as they have the simplest derivatives among all nonlinear functions. 
The ReLU function is proved to be useful in training deep neural network architectures
The Leaky ReLU function can resolve the vanishing gradient problem 
which can accelerate the training in some occasions. 
The optimizers Adamax and Proximal Adagrad are stochastic gradient descent (SGD) 
based methods commonly used in neural network training \cite{adam}. 
In both experiments, we trained our network using Python API 
Tensorflow and Keras \cite{chollet2015keras}. 

Once a neural network is built on training, 
it can be used to predict the output given a new input. 
The accuracy of the predictions is essential in making the network useful. 
In our experiments, we use of $M$ sample pairs, 
which are not used in training the network, 
to examine the predictive power of our network. 
On these sample pairs, referred to as the testing set, 
we compare the prediction and the exact output 
and compute the mismatch in some suitable metric. 
Here, we summarize the metric used in our numerical experiment. 
For the multiscale basis functions, we compute the relative error in $L^2$ and $H^1$ norm, i.e. 
\begin{equation}
\begin{split}
e_{L^2}(\kappa_{N+j}) & =  \left(
\dfrac{\int_{\Omega} \left\vert \phi_m^{\omega_i}(\kappa_{N+j}) - 
\phi_m^{\omega_i,\text{pred}}(\kappa_{N+j}) \right\vert^2}
{ \int_{\Omega} \left\vert \phi_m^{\omega_i}(\kappa_{N+j}) 
\right\vert^2}\right)^\frac{1}{2}, \\
e_{H^1}(\kappa_{N+j}) & = \left(
\dfrac{\int_{\Omega} \left\vert \nabla \phi_m^{\omega_i}(\kappa_{N+j}) - 
\nabla \phi_m^{\omega_i,\text{pred}}(\kappa_{N+j}) \right\vert^2}
{ \int_{\Omega} \left\vert \nabla \phi_m^{\omega_i}(\kappa_{N+j}) 
\right\vert^2}\right)^\frac{1}{2}.
\end{split}
\end{equation}
For the local stiffness matrices, we compute the relative error in entrywise $\ell^2$, 
entrywise $\ell^\infty$ and Frobenius norm, i.e. 
\begin{equation}
\begin{split}
e_{\ell^2}(\kappa_{N+j}) & = 
\dfrac{\| A_c^{K_l}(\kappa_{N+j}) - A_c^{K_l,\text{pred}}(\kappa_{N+j}) \|_2}
{\| A_c^{K_l}(\kappa_{N+j})\|_2}, \\
e_{\ell^\infty}(\kappa_{N+j}) & = 
\dfrac{\| A_c^{K_l}(\kappa_{N+j}) - A_c^{K_l,\text{pred}}(\kappa_{N+j}) \|_\infty}
{\| A_c^{K_l}(\kappa_{N+j})\|_\infty}, \\
e_{F}(\kappa_{N+j}) & = 
\dfrac{\| A_c^{K_l}(\kappa_{N+j}) - A_c^{K_l,\text{pred}}(\kappa_{N+j}) \|_F}
{\| A_c^{K_l}(\kappa_{N+j})\|_F}. 
\end{split}
\end{equation}
A more important measure of the usefulness of the trained neural network is the 
predicted multiscale solution $u_{ms}^{\text{pred}}(\kappa)$ 
given by \eqref{eq:coarse_system_pred}--\eqref{eq:downscale_pred}. 
We compare the predicted solution to $u_{ms}$ defined by 
\eqref{eq:coarse_system}--\eqref{eq:downscale}, 
and compute the relative error in $L^2$ and energy norm, i.e. 
\begin{equation}
\begin{split}
e_{L^2}(\kappa_{N+j}) & = \left(
\dfrac{\int_{\Omega} \left\vert u_{ms}(\kappa_{N+j}) - 
u_{ms}^{\text{pred}}(\kappa_{N+j}) \right\vert^2}
{ \int_{\Omega} \left\vert u_{ms}(\kappa_{N+j}) 
\right\vert^2}\right)^\frac{1}{2}, \\
e_{a}(\kappa_{N+j}) & = \left(
\dfrac{\int_{\Omega} \kappa_j \left\vert \nabla u_{ms}(\kappa_{N+j}) - 
\nabla u_{ms}^{\text{pred}}(\kappa_{N+j}) \right\vert^2}
{ \int_{\Omega}  \kappa_j \left\vert \nabla u_{ms}(\kappa_{N+j}) 
\right\vert^2}\right)^\frac{1}{2}.
\end{split}
\end{equation}

\subsection{Experiment 1}

\begin{figure}[h!]
	\centering
	\begin{subfigure}{.4\textwidth}
		\centering
		\includegraphics[width=.8\linewidth]{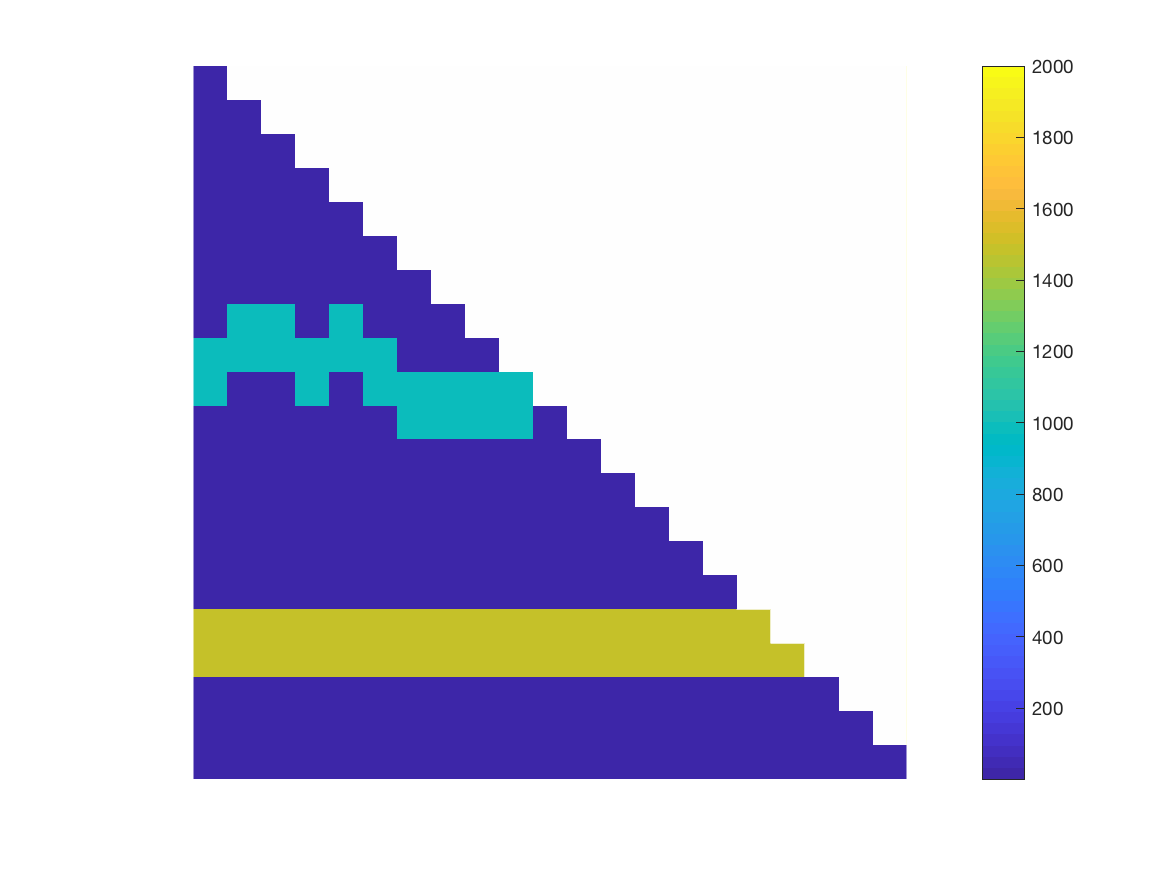}
		\label{fig:kappa_one_4}
	\end{subfigure}
	\begin{subfigure}{.4\textwidth}
		\centering
		\includegraphics[width=.8\linewidth]{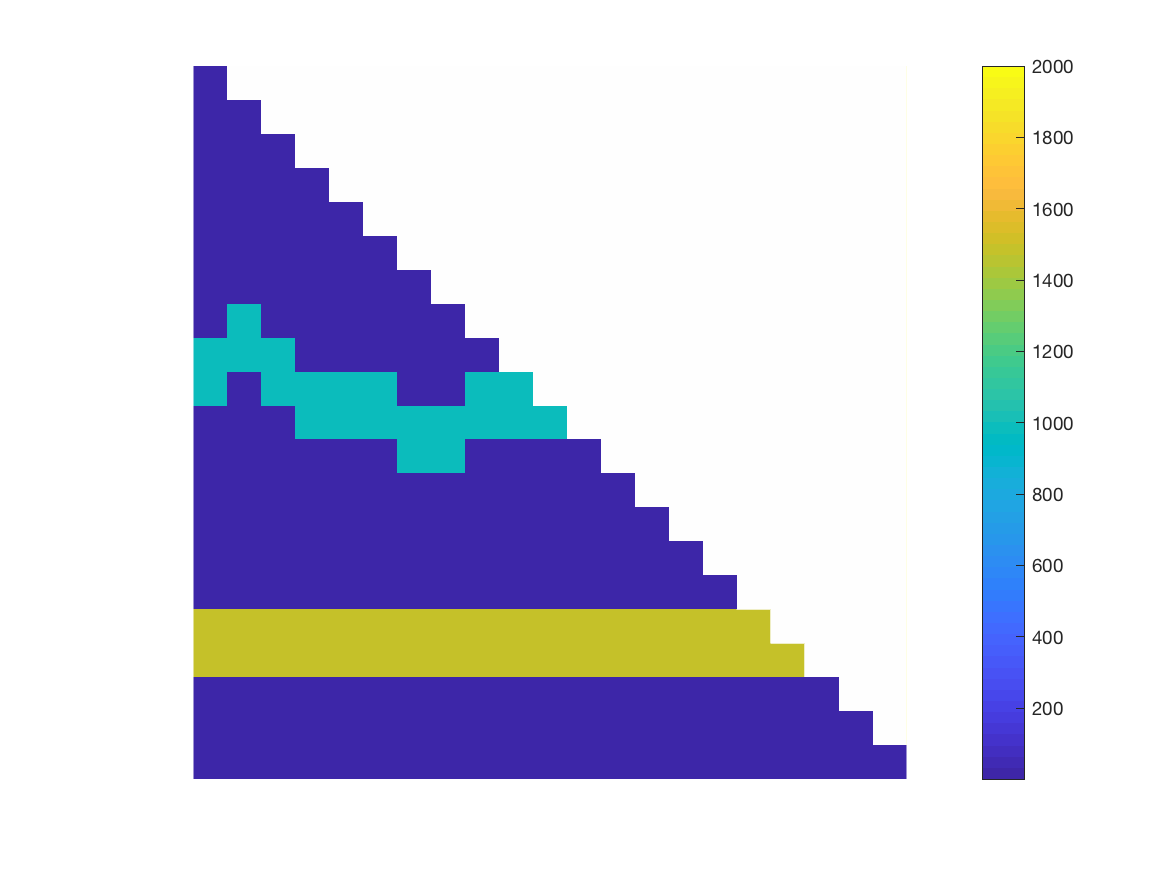}
		\label{fig:kappa_one_991}
	\end{subfigure}
	
	\begin{subfigure}{.4\textwidth}
		\centering
		\includegraphics[width=.8\linewidth]{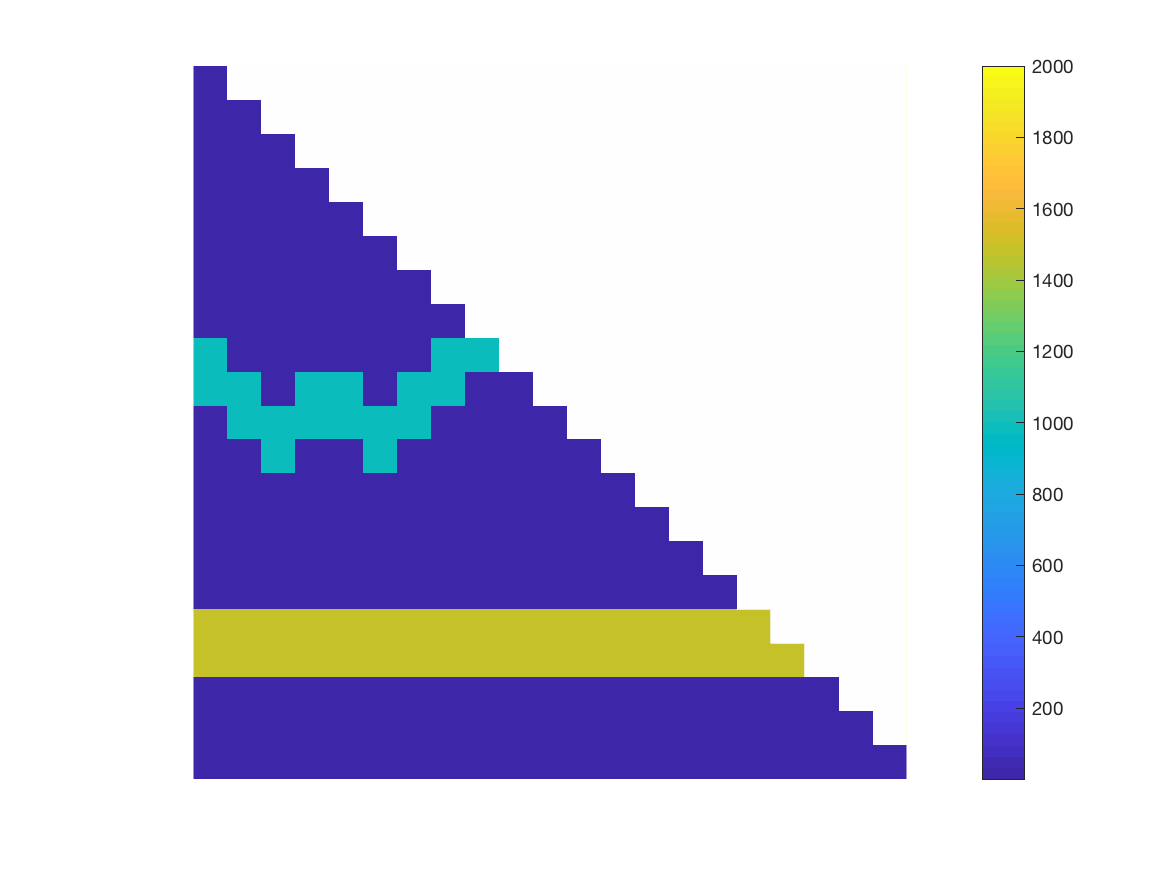}
		\label{fig:kappa_one_1443}
	\end{subfigure}
	\begin{subfigure}{.4\textwidth}
		\centering
		\includegraphics[width=.8\linewidth]{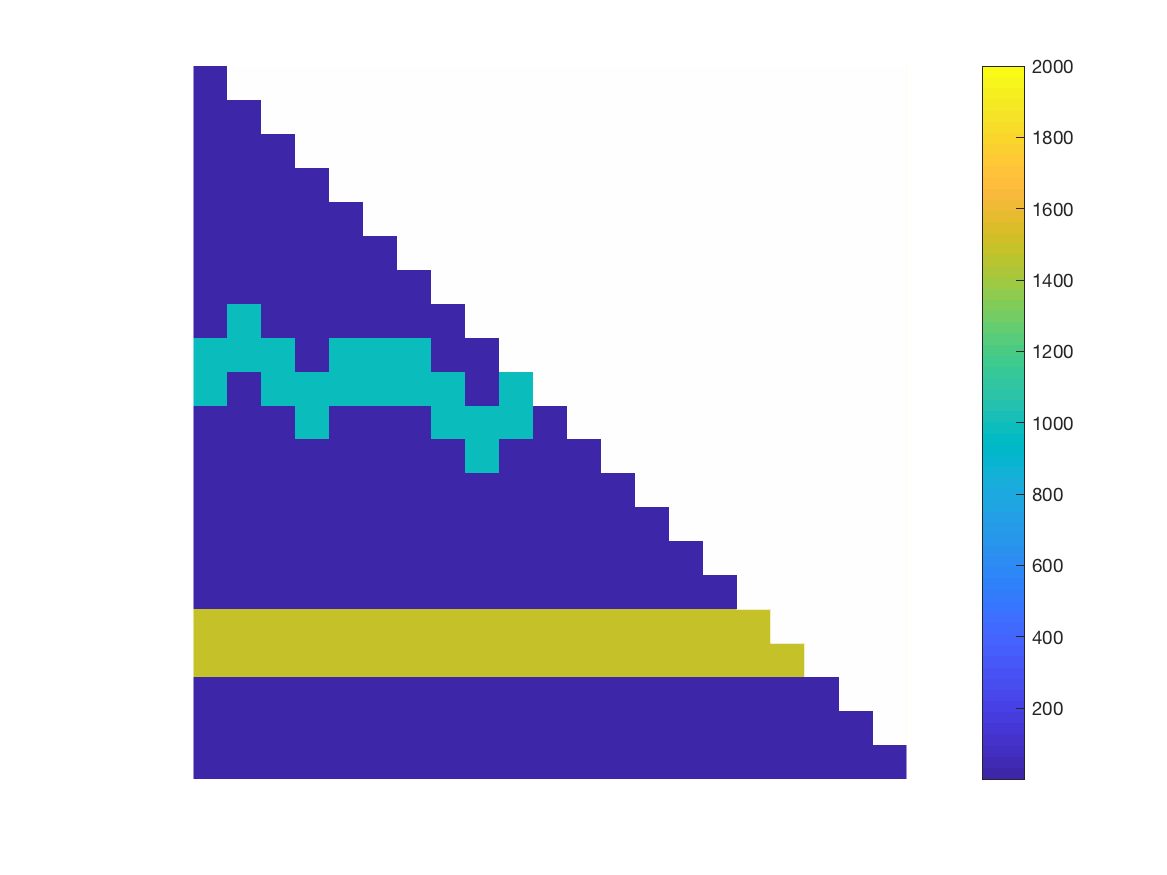}
		\label{fig:kappa_one_1528}
	\end{subfigure}
	
	\begin{subfigure}{.4\textwidth}
		\centering
		\includegraphics[width=.8\linewidth]{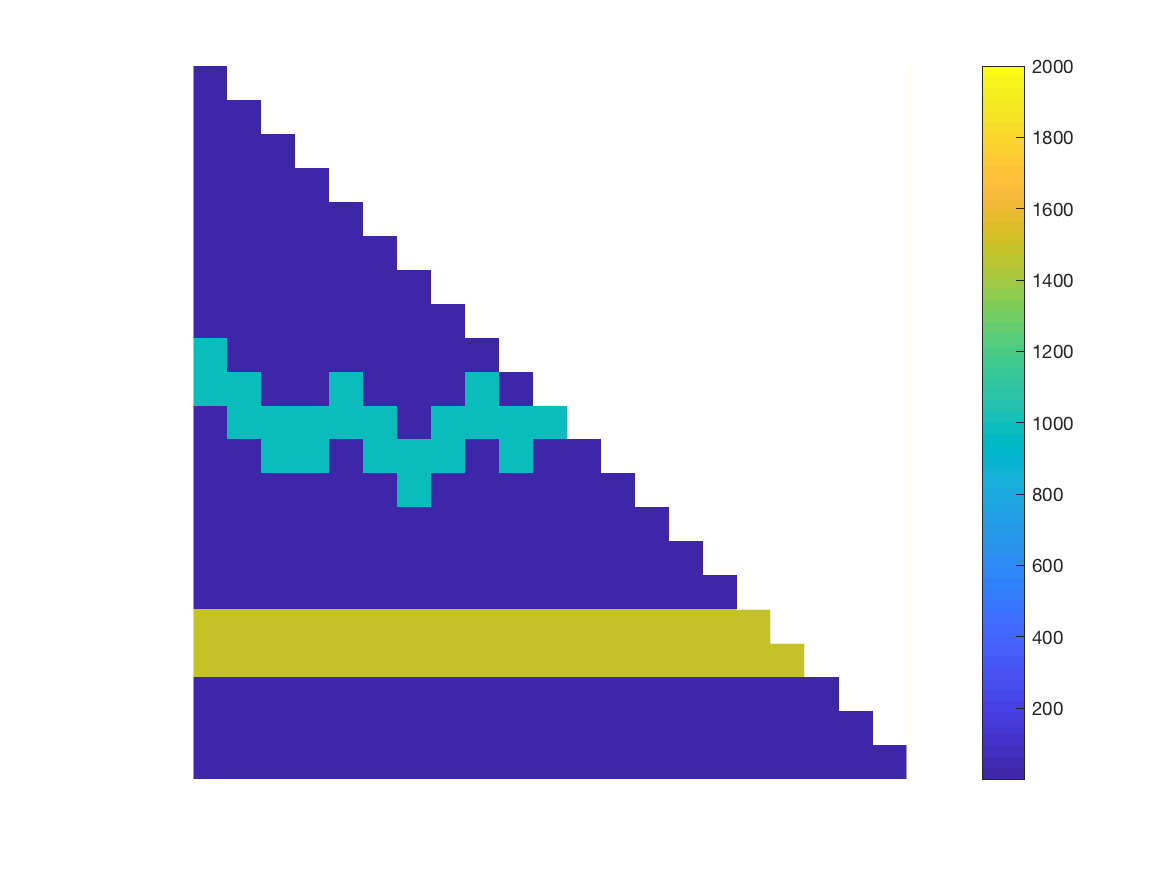}
		\label{fig:kappa_one_1803}
	\end{subfigure}
	\begin{subfigure}{.4\textwidth}
		\centering
		\includegraphics[width=.8\linewidth]{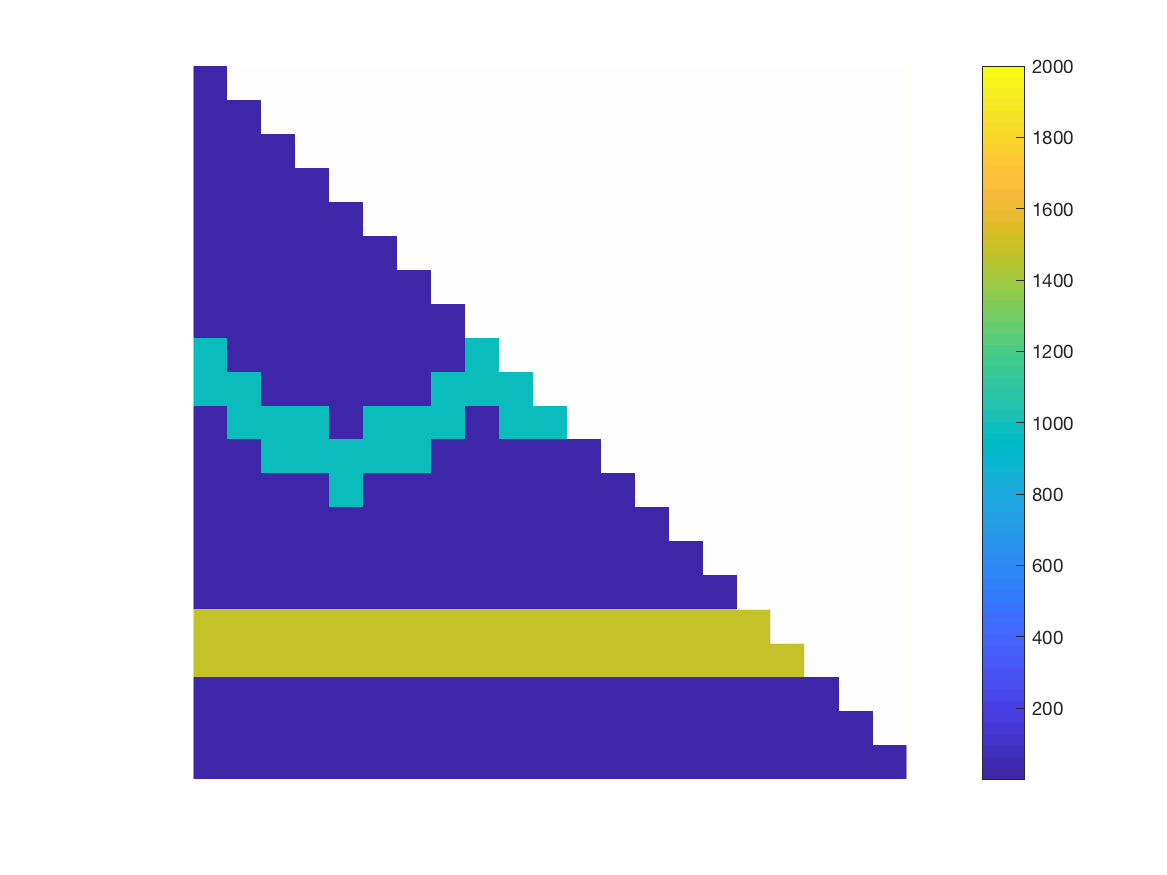}
		\label{fig:kappa_one_1891}
	\end{subfigure}
	\caption{Samples of permeability fields in the target block $K_0$ in Experiment 1.}
	\label{fig:kappa_one_samples}
\end{figure}

In this experiment, we consider curved channelized 
permeability fields. Each permeability field contains 
a straight channel and a curved channel. 
The straight channel is fixed and 
the curved channel strikes the boundary of the target cell $K_0$ 
at the same points. The curvature of the sine-shaped channel inside $K_0$
varies among all realizations. 
We generate 2000 realizations of 
permeability fields, where the permeability coefficients are fixed. 
Samples of permeability fields are depicted in Figure~\ref{fig:kappa_one_samples}. 
Among the 2000 realizations, 
1980 sample pairs are randomly chosen and used as training samples, 
and the remaining 20 sample pairs are used as testing samples.

For each realization, we compute the 
local multiscale basis functions and local coarse-scale stiffness matrix. 
In building the local snapshot space, we solve for harmonic extension of 
all the fine-grid boundary conditions.
Local multiscale basis functions are then constructed by solving the spectral problem 
and multiplied the spectral basis functions with the multiscale partition of unity functions. 
With the offline space constructed, 
we compute the coarse-scale stiffness matrix. 
We use the training samples to build deep neural networks 
for approximating these GMsFEM quantities, 
and examine the performance of the approximations 
on the testing set. 

Tables~\ref{tab:exp1b1e}--\ref{tab:exp1me} record 
the error of the prediction by the neural networks in each testing sample 
and the mean error 
measured in the defined metric. 
It can be seen that the prediction are of high accuracy. 
This is vital in ensuring the predicted GMsFEM solver useful. 
Table~\ref{tab:exp1se} records the error of the multiscale solution 
in each testing sample and the mean error
using our proposed method. 
It can be observed that using the predicted GMsFEM solver, 
we obtain a good approximation of the multiscale solution 
compared with the exact GMsFEM solver.

\begin{table}[ht!]
\centering
\begin{tabular}{c|cc|cc|cc}
Sample & \multicolumn{2}{|c|}{$\omega_1$} & \multicolumn{2}{|c|}{$\omega_2$} & \multicolumn{2}{|c}{$\omega_3$} \\
$j$ & $e_{L^2}$ & $e_{H^1}$ & $e_{L^2}$ & $e_{H^1}$ & $e_{L^2}$ & $e_{H^1}$ \\
\hline
1 & 0.47\% & 3.2\% & 0.40\% & 3.6\% & 0.84\% & 5.1\% \\ 
2 & 0.45\% & 4.4\% & 0.39\% & 3.3\% & 1.00\% & 6.3\% \\ 
3 & 0.34\% & 2.3\% & 0.40\% & 3.1\% & 0.88\% & 4.3\% \\ 
4 & 0.35\% & 4.2\% & 0.43\% & 5.4\% & 0.94\% & 6.6\% \\ 
5 & 0.35\% & 3.3\% & 0.37\% & 3.9\% & 0.90\% & 6.1\% \\ 
6 & 0.51\% & 4.7\% & 0.92\% & 12.0\% & 2.60\% & 19.0\% \\ 
7 & 0.45\% & 4.1\% & 0.38\% & 3.2\% & 1.00\% & 6.4\% \\ 
8 & 0.31\% & 3.4\% & 0.43\% & 5.5\% & 1.10\% & 7.7\% \\ 
9 & 0.25\% & 2.2\% & 0.46\% & 5.6\% & 1.10\% & 6.2\% \\ 
10 & 0.31\% & 3.5\% & 0.42\% & 4.5\% & 1.30\% & 7.6\% \\ 
\hline
Mean & 0.38\% & 3.5\% & 0.46\% & 5.0\% & 1.17\% & 7.5\% \\ 
\end{tabular}
\caption{Percentage error of multiscale basis functions $\phi_1^{\omega_i}$ in Experiment 1.}
\label{tab:exp1b1e}
\end{table}

\begin{table}[ht!]
\centering
\begin{tabular}{c|cc|cc|cc}
Sample & \multicolumn{2}{|c|}{$\omega_1$} & \multicolumn{2}{|c|}{$\omega_2$} & \multicolumn{2}{|c}{$\omega_3$} \\
$j$ & $e_{L^2}$ & $e_{H^1}$ & $e_{L^2}$ & $e_{H^1}$ & $e_{L^2}$ & $e_{H^1}$ \\
\hline
1 & 0.47\% & 4.2\% & 0.40\% & 1.4\% & 0.32\% & 1.1\% \\ 
2 & 0.57\% & 3.2\% & 0.31\% & 1.4\% & 0.30\% & 1.1\% \\ 
3 & 0.58\% & 2.7\% & 0.31\% & 1.4\% & 0.33\% & 1.1\% \\ 
4 & 0.59\% & 3.6\% & 0.13\% & 1.3\% & 0.32\% & 1.1\% \\ 
5 & 0.53\% & 4.0\% & 0.51\% & 1.6\% & 0.27\% & 1.0\% \\ 
6 & 0.85\% & 4.3\% & 0.51\% & 2.1\% & 0.29\% & 1.3\% \\ 
7 & 0.50\% & 2.7\% & 0.22\% & 1.5\% & 0.29\% & 1.0\% \\ 
8 & 0.43\% & 4.5\% & 0.61\% & 1.9\% & 0.35\% & 1.1\% \\ 
9 & 0.71\% & 2.9\% & 0.14\% & 1.4\% & 0.27\% & 1.1\% \\ 
10 & 0.66\% & 4.4\% & 0.53\% & 1.8\% & 0.26\% & 1.1\% \\ 
\hline
Mean & 0.59\% & 3.6\% & 0.37\% & 1.6\% & 0.30\% & 1.1\% \\ 
\end{tabular}
\caption{Percentage error of multiscale basis functions $\phi_2^{\omega_i}$ in Experiment 1.}
\label{tab:exp1b2e}
\end{table}

\begin{table}[ht!]
\centering
\begin{tabular}{c|cc}
Sample $j$ & $e_{\ell^2}$     & $e_{F}$\\
\hline
1& 	0.67\%	& 0.84\% \\
2& 	0.37\%	& 0.37\%\\
3& 0.32\%	& 0.38\%\\
4& 	1.32\%	& 1.29\%\\
5& 0.51\%	& 0.59\%\\
6& 4.43\%	& 4.28\%\\
7& 	0.34\%	& 0.38\%\\
8& 0.86\%	& 1.04\%\\
9& 	1.00\%	& 0.97\%\\
10& 	0.90\%	& 1.08\%\\
\hline
Mean & 0.76\% & 0.81\%
\end{tabular}
\caption{Percentage error of the local stiffness matrix $A_c^{K_0}$ in Experiment 1.}
\label{tab:exp1me}%
\end{table}

\begin{table}[ht!]
\centering
\begin{tabular}{c|cc}
Sample $j$ & $e_{L^2}$  & $e_{a}$\\
\hline
1 & 0.31\% & 4.58\% \\ 
2 & 0.30\% & 4.60\% \\ 
3 & 0.30\% & 4.51\% \\ 
4 & 0.27\% & 4.60\% \\ 
5 & 0.29\% & 4.56\% \\ 
6 & 0.47\% & 4.67\% \\ 
7 & 0.39\% & 4.70\% \\ 
8 & 0.30\% & 4.63\% \\ 
9 & 0.35\% & 4.65\% \\ 
10 & 0.31\% & 4.65\% \\ 
\hline
Mean & 0.33\% & 4.62\% \\ 
\end{tabular}
\caption{Percentage error of multiscale solution $u_{ms}$ in Experiment 1.}
\label{tab:exp1se}%
\end{table}

\subsection{Experiment 2}\label{section:sine_samples}

\begin{figure}[ht!]
	\centering
	\begin{subfigure}{.32\textwidth}
		\centering
		\includegraphics[width=.8\linewidth]{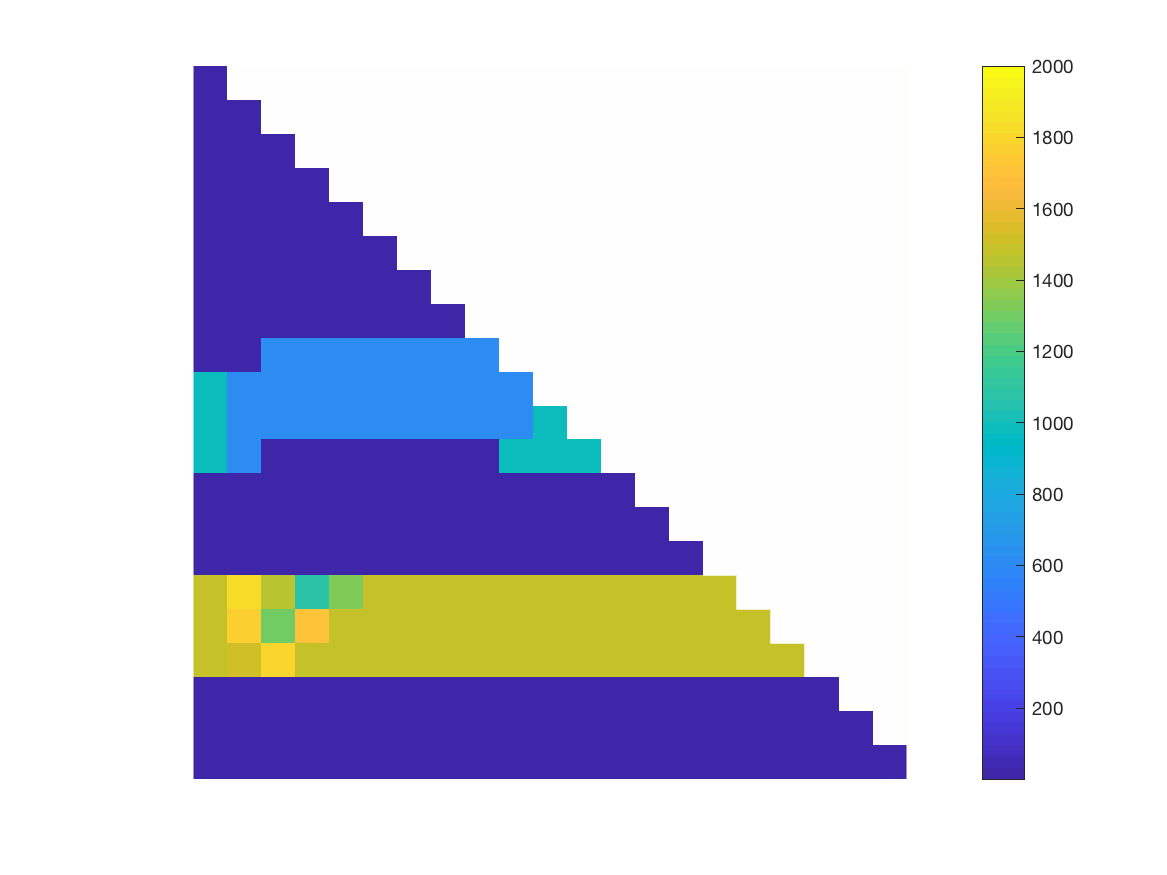}
		\label{fig:sine_sample5}
	\end{subfigure}%
	\begin{subfigure}{.32\textwidth}
		\centering
		\includegraphics[width=.8\linewidth]{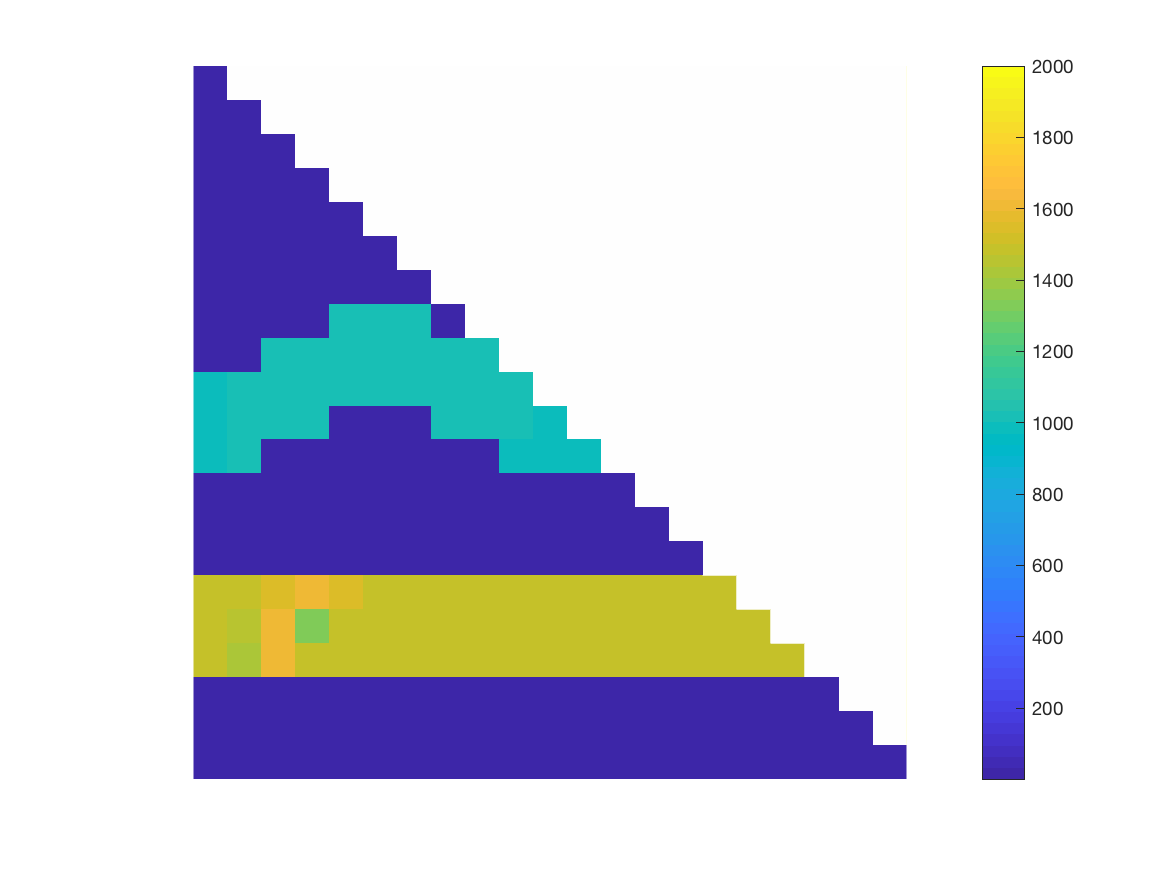}
		\label{fig:sine_sample541}
	\end{subfigure}
	\begin{subfigure}{.32\textwidth}
		\centering
		\includegraphics[width=.8\linewidth]{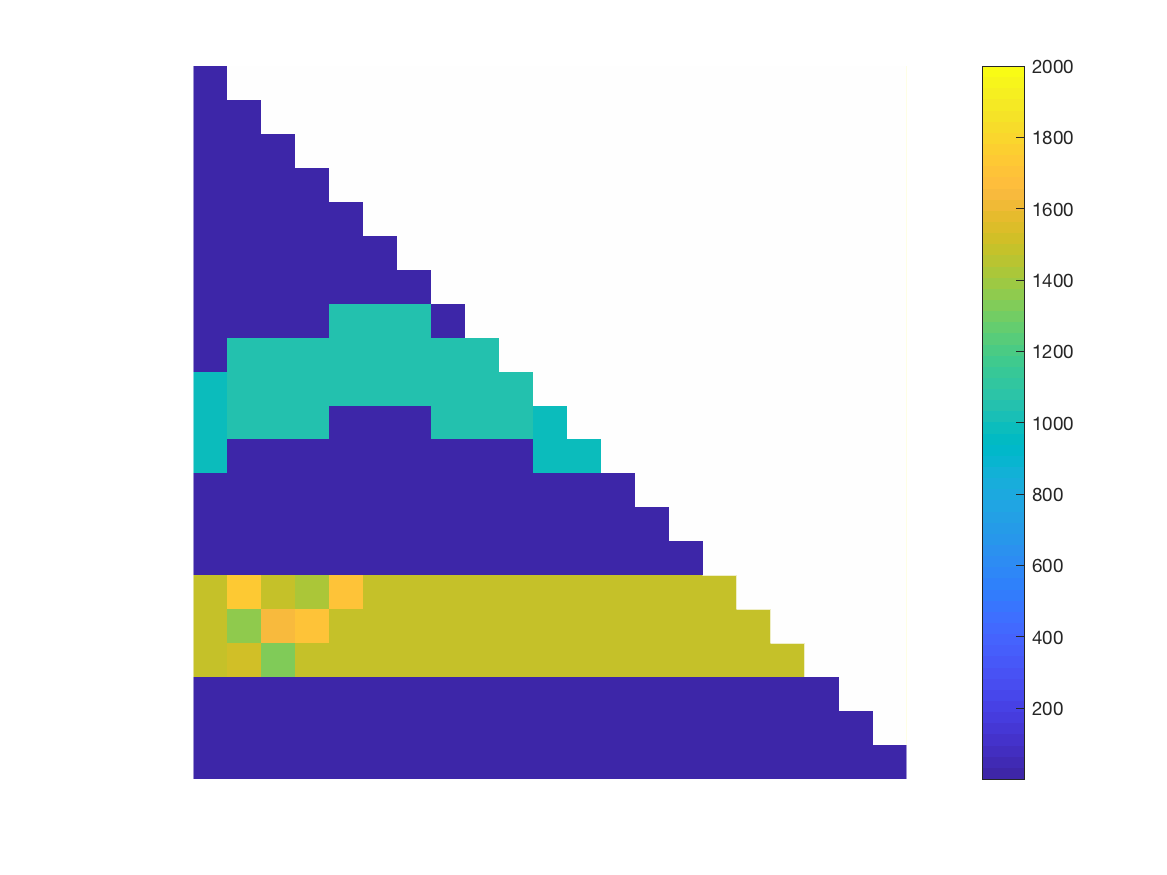}
		\label{fig:sine_sample1113}
	\end{subfigure} \\
	\begin{subfigure}{.32\textwidth}
		\centering
		\includegraphics[width=.8\linewidth]{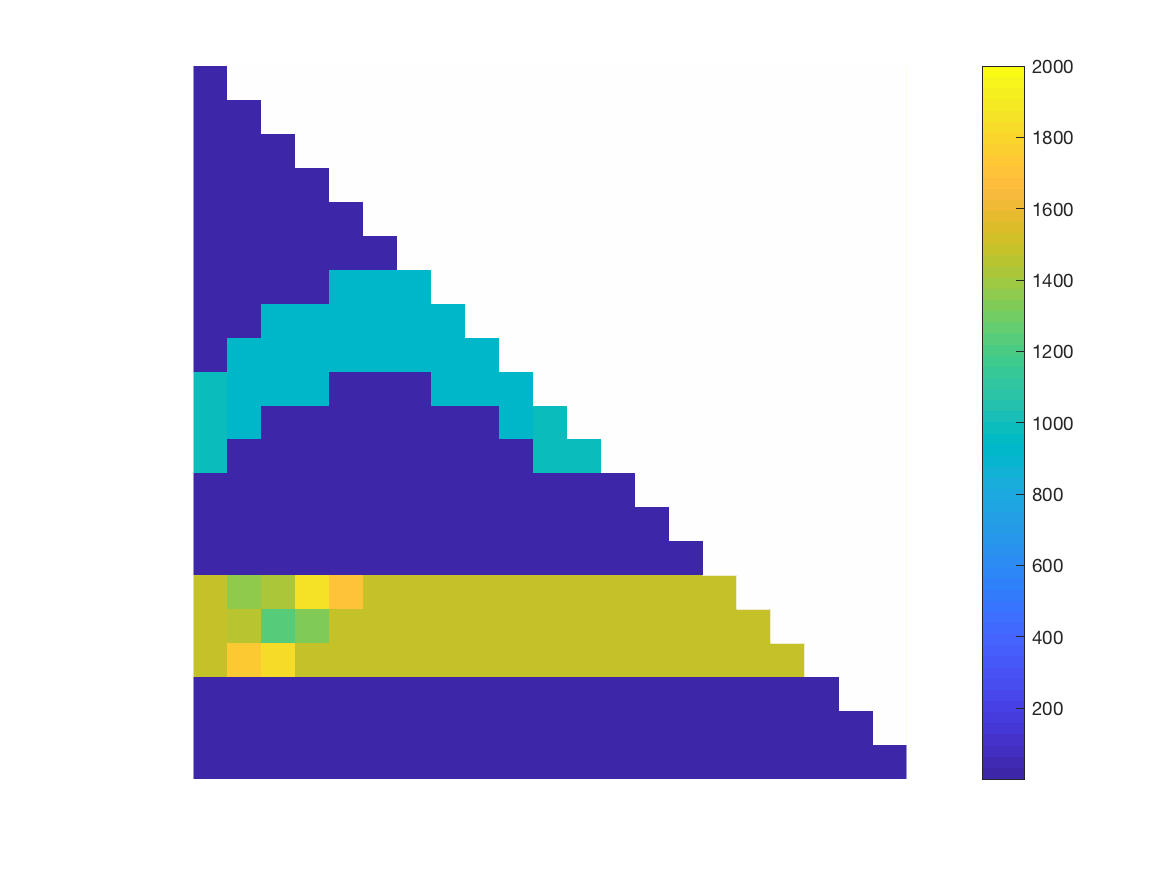}
		\label{fig:sine_sample1906}
	\end{subfigure}
	\begin{subfigure}{.32\textwidth}
		\centering
		\includegraphics[width=.8\linewidth]{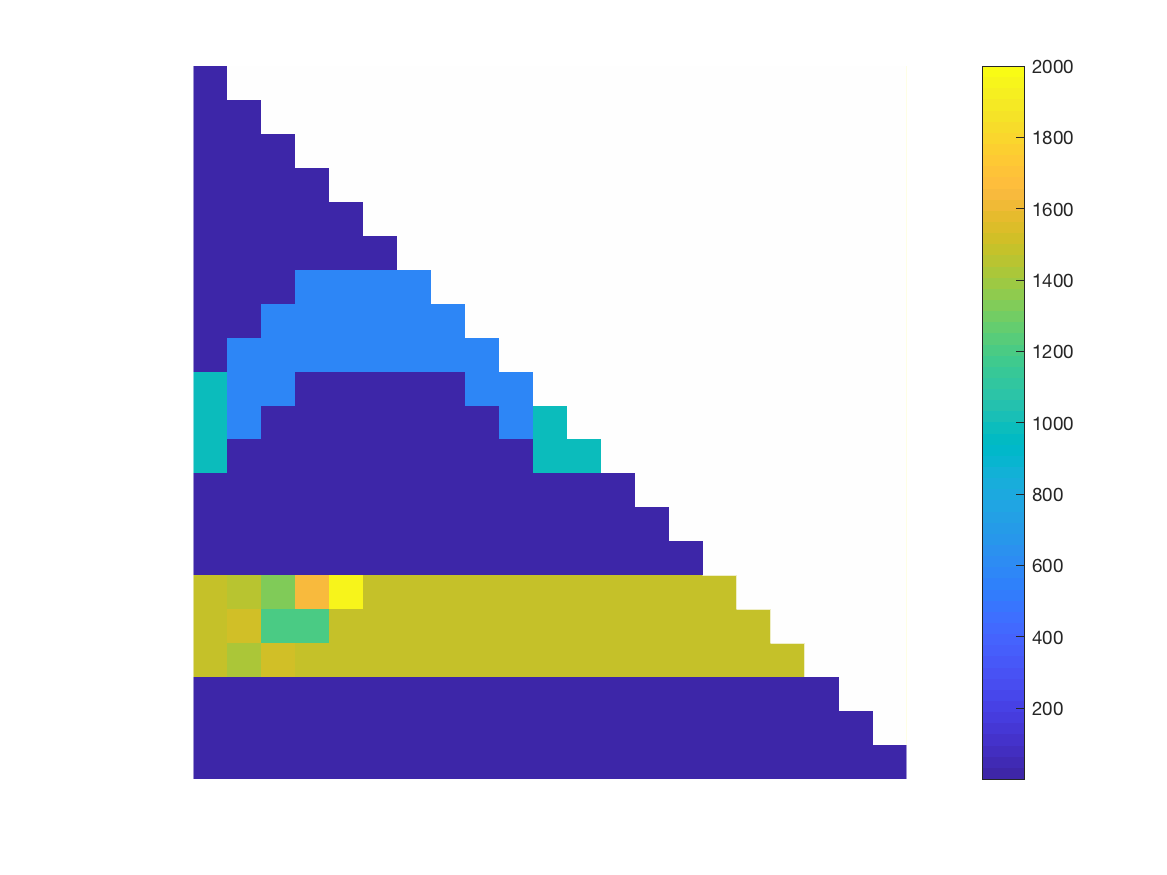}
		\label{fig:sine_sample2254}
	\end{subfigure}
	
	\caption{Samples of permeability fields in the target block $K_0$ in Experiment 2.}
	\label{fig:sine_samples}
\end{figure}

In this experiment, we consider sine-shaped channelized 
permeability fields. Each permeability field contains 
a straight channel and a sine-shaped channel. 
There are altogether 5 channel configurations, 
where the straight channel is fixed and 
the sine-shaped channel strikes the boundary of the target cell $K_0$ 
at the same points. The curvature of the sine-shaped channel inside $K_0$
varies among these configurations. 
For each channel configuration, we generate 500 realizations of 
permeability fields, where the permeability coefficients follow random distributions. 
Samples of permeability fields are depicted in Figure~\ref{fig:sine_samples}. 
Among the 2500 realizations, 
2475 sample pairs are randomly chosen and used as training samples, 
and the remaining 25 sample pairs are used as testing samples. 

Next, for each realization, we compute the 
local multiscale basis functions and local coarse-scale stiffness matrix. 
In building the local snapshot space, we solve for harmonic extension of 
randomized fine-grid boundary conditions, 
so as to reduce the number of local problems to be solved. 
Local multiscale basis functions are then constructed by solving the spectral problem 
and multiplied the spectral basis functions with the multiscale partition of unity functions. 
With the offline space constructed, 
we compute the coarse-scale stiffness matrix. 
We use the training samples to build deep neural networks 
for approximating these GMsFEM quantities, 
and examine the performance of the approximations 
on the testing set. 

\begin{figure}[ht!]
	\begin{subfigure}{.3\textwidth}
		\centering
		\includegraphics[width=.8\linewidth]{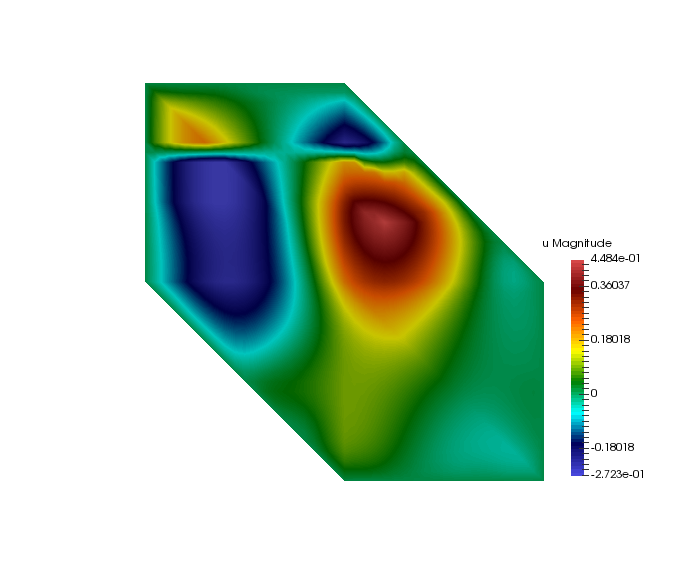}
		\caption{Exact $\phi_1^{\omega_1}$}
	\end{subfigure}
	\begin{subfigure}{.3\textwidth}
		\centering
		\includegraphics[width=.8\linewidth]{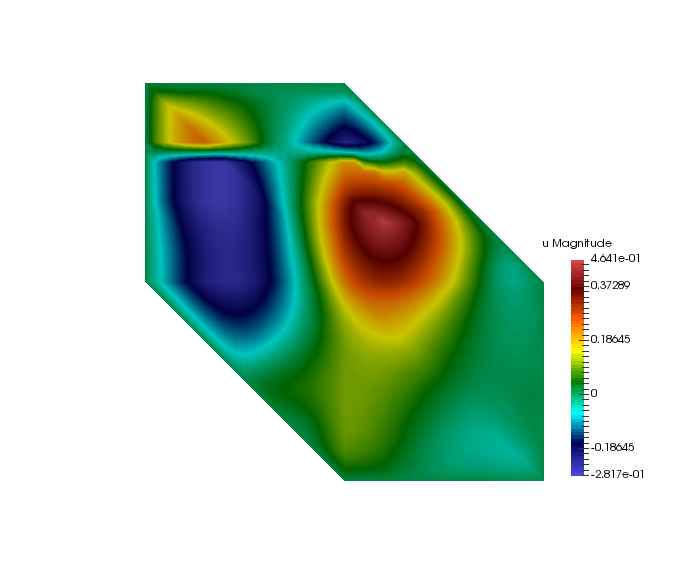}
		\caption{Prediction $\phi_1^{\omega_1,\text{pred}}$}
	\end{subfigure}
	\begin{subfigure}{.3\textwidth}
		\centering
		\includegraphics[width=.8\linewidth]{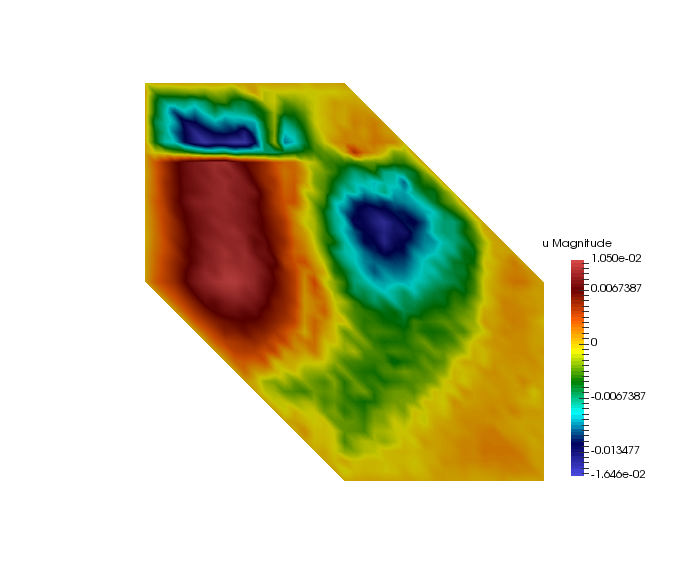}
		\caption{Difference $\phi_1^{\omega_1} - \phi_1^{\omega_1,\text{pred}}$}
	\end{subfigure}
	
	\begin{subfigure}{.3\textwidth}
		\centering
		\includegraphics[width=.8\linewidth]{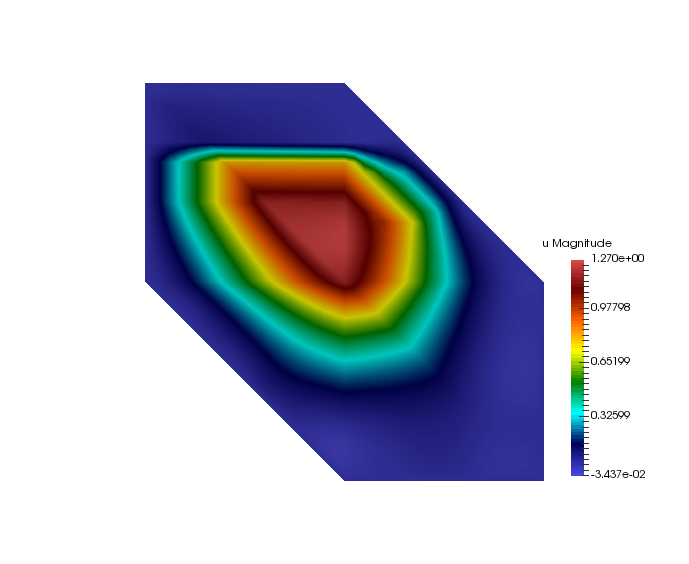}
		\caption{Exact $\phi_2^{\omega_1}$}
	\end{subfigure}
	\begin{subfigure}{.3\textwidth}
		\centering
		\includegraphics[width=.8\linewidth]{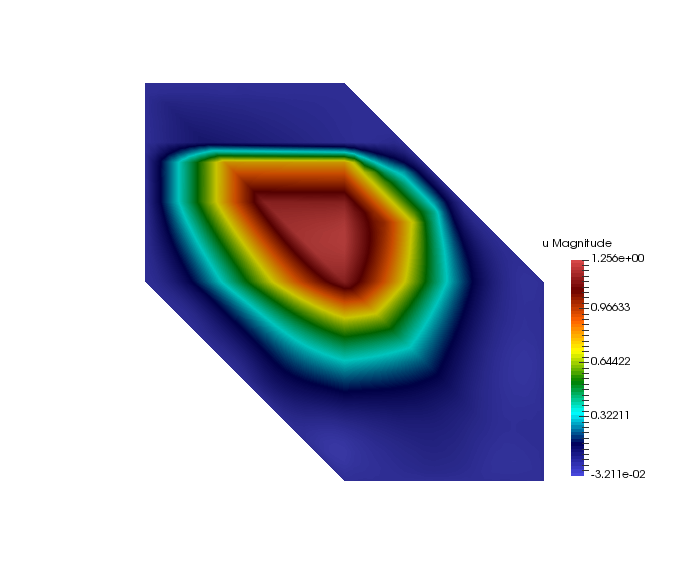}
		\caption{Prediction $\phi_2^{\omega_1,\text{pred}}$}
	\end{subfigure}
	\begin{subfigure}{.3\textwidth}
		\centering
		\includegraphics[width=.8\linewidth]{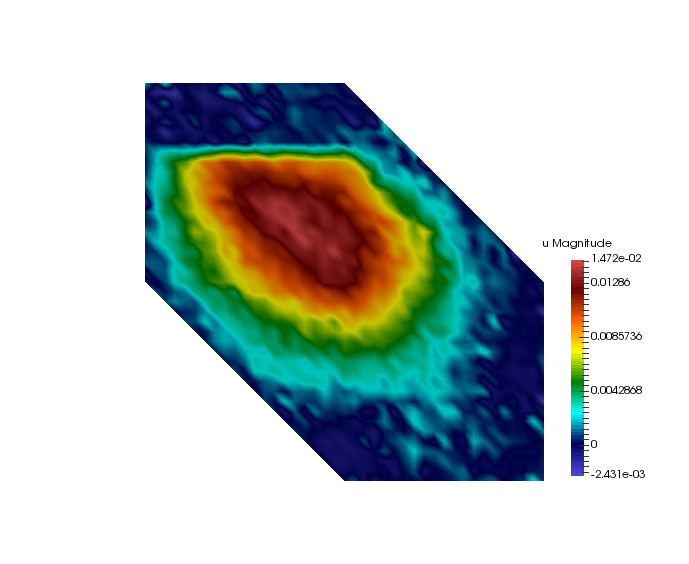}
		\caption{Difference $\phi_2^{\omega_1} - \phi_2^{\omega_1,\text{pred}}$}
	\end{subfigure}
	\caption{Exact multiscale basis functions $\phi_m^{\omega_1}$ and predicted multiscale basis functions $\phi_m^{\omega_1,\text{pred}}$ in the coarse neighborhood $\omega_1$ in Experiment 2.}
	\label{fig:88_basis_rand}
\end{figure}

\begin{figure}[ht!]
	\begin{subfigure}{.3\textwidth}
		\centering
		\includegraphics[width=.8\linewidth]{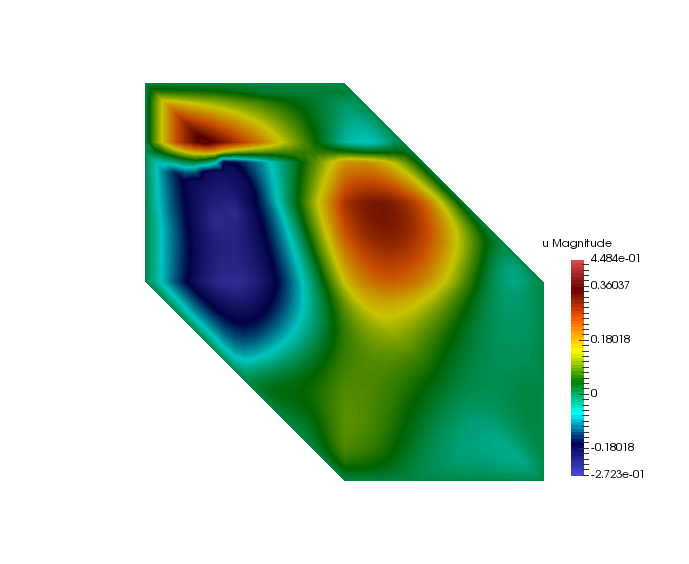}
		\caption{Exact $\phi_1^{\omega_2}$}
	\end{subfigure}
	\begin{subfigure}{.3\textwidth}
		\centering
		\includegraphics[width=.8\linewidth]{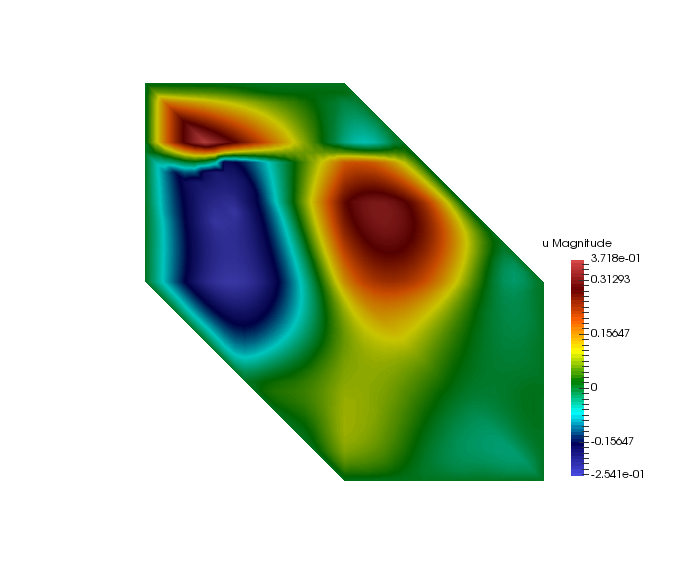}
		\caption{Prediction $\phi_1^{\omega_2,\text{pred}}$}
	\end{subfigure}
	\begin{subfigure}{.3\textwidth}
		\centering
		\includegraphics[width=.8\linewidth]{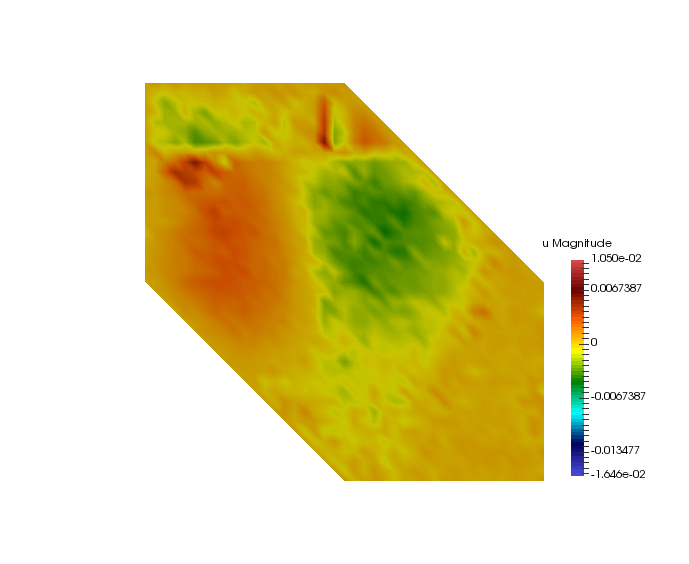}
		\caption{Difference $\phi_1^{\omega_2} - \phi_1^{\omega_2,\text{pred}}$}
	\end{subfigure}
	
	\begin{subfigure}{.3\textwidth}
		\centering
		\includegraphics[width=.8\linewidth]{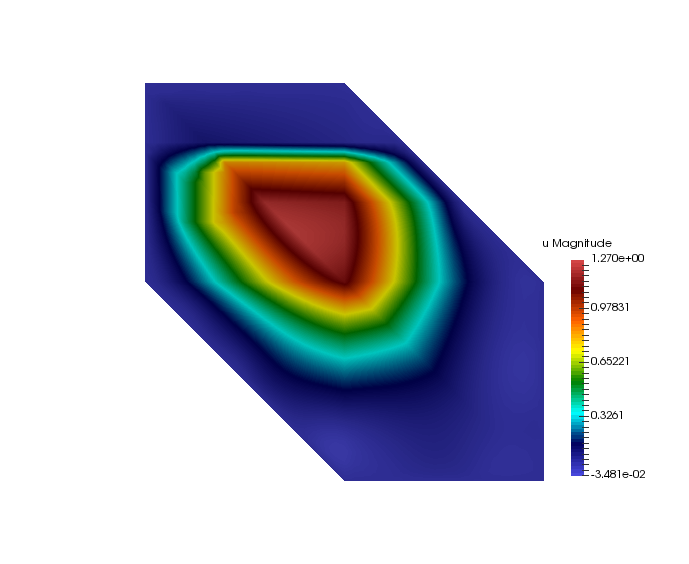}
		\caption{Exact $\phi_2^{\omega_2}$}
	\end{subfigure}
	\begin{subfigure}{.3\textwidth}
		\centering
		\includegraphics[width=.8\linewidth]{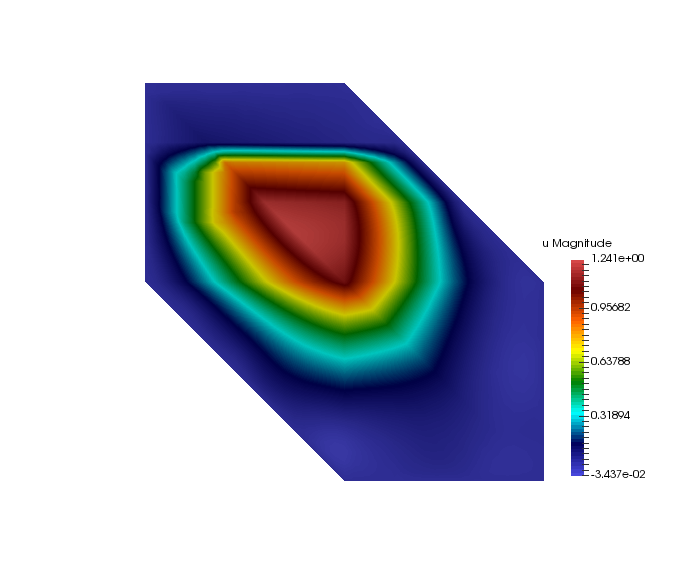}
		\caption{Prediction $\phi_2^{\omega_2,\text{pred}}$}
	\end{subfigure}
	\begin{subfigure}{.3\textwidth}
		\centering
		\includegraphics[width=.8\linewidth]{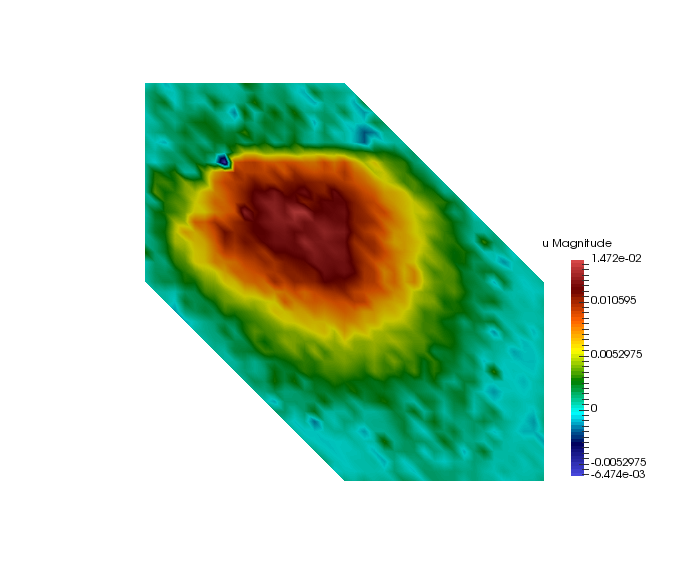}
		\caption{Difference $\phi_2^{\omega_2} - \phi_2^{\omega_2,\text{pred}}$}
	\end{subfigure}
	\caption{Exact multiscale basis functions $\phi_m^{\omega_2}$ and predicted multiscale basis functions $\phi_m^{\omega_2,\text{pred}}$ in the coarse neighborhood $\omega_2$ in Experiment 2.}
	\label{fig:89_basis_rand}
\end{figure}

\begin{figure}[ht!]
	\begin{subfigure}{.3\textwidth}
		\centering
		\includegraphics[width=.8\linewidth]{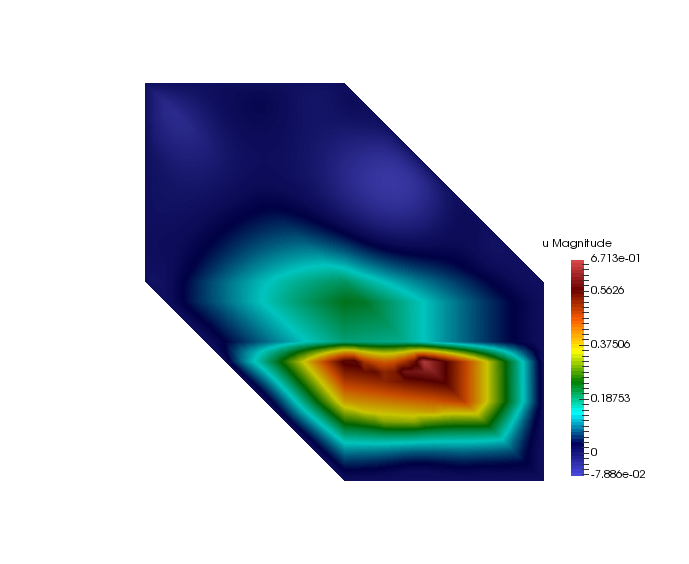}
		\caption{Exact $\phi_1^{\omega_3}$}
	\end{subfigure}
	\begin{subfigure}{.3\textwidth}
		\centering
		\includegraphics[width=.8\linewidth]{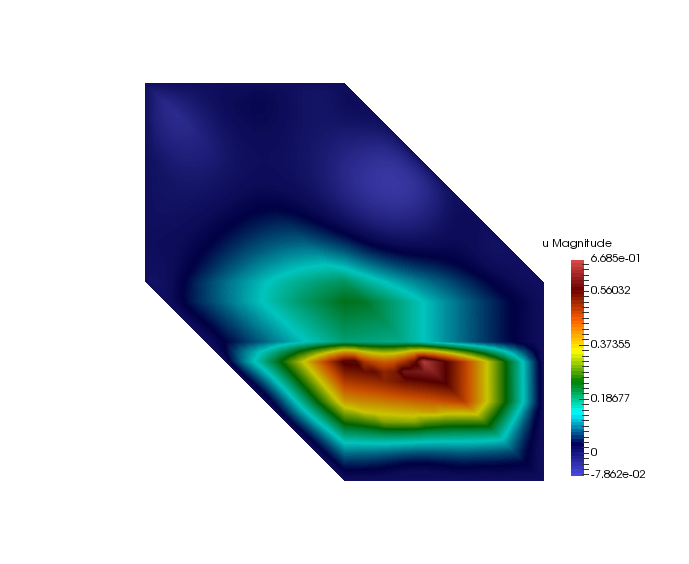}
		\caption{Prediction $\phi_1^{\omega_3,\text{pred}}$}
	\end{subfigure}
	\begin{subfigure}{.3\textwidth}
		\centering
		\includegraphics[width=.8\linewidth]{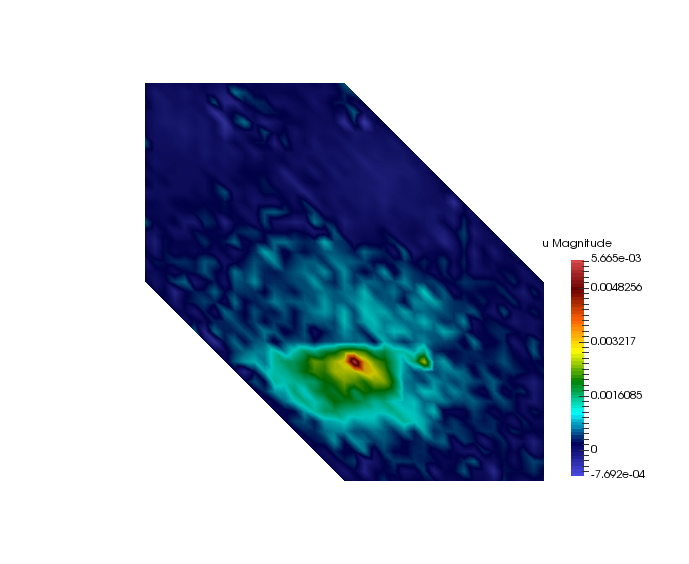}
		\caption{Difference $\phi_1^{\omega_3} - \phi_1^{\omega_3,\text{pred}}$}
	\end{subfigure}
	
	\begin{subfigure}{.3\textwidth}
		\centering
		\includegraphics[width=.8\linewidth]{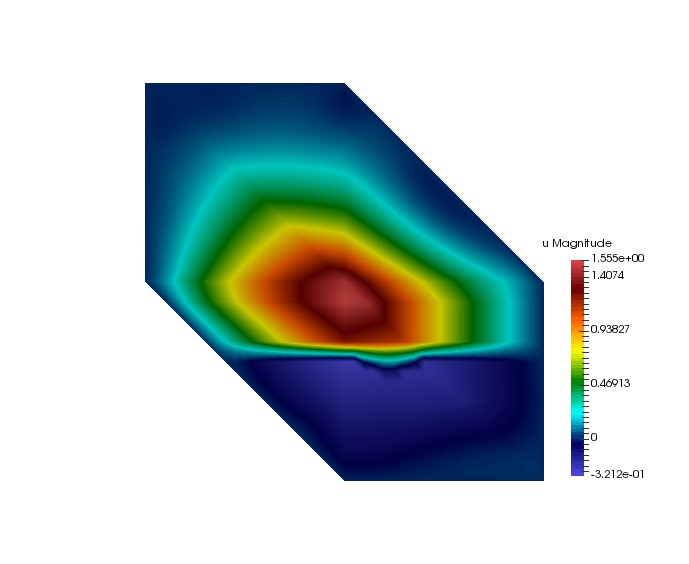}
		\caption{Exact $\phi_2^{\omega_3}$}
	\end{subfigure}
	\begin{subfigure}{.3\textwidth}
		\centering
		\includegraphics[width=.8\linewidth]{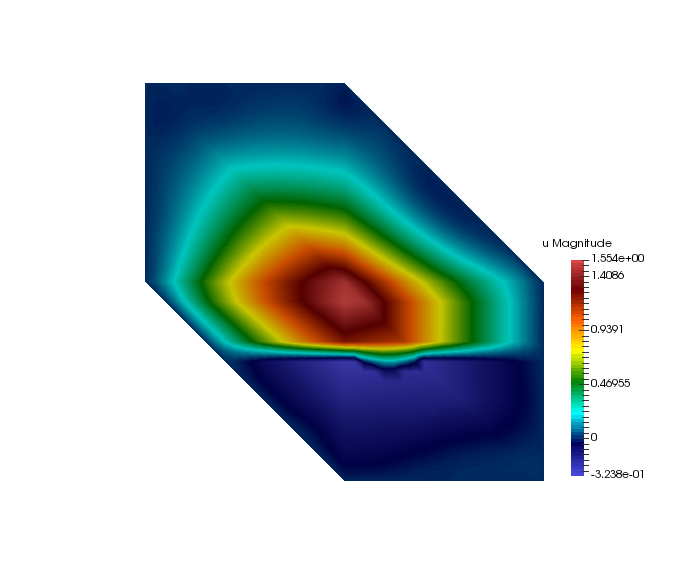}
		\caption{Prediction $\phi_2^{\omega_3,\text{pred}}$}
	\end{subfigure}
	\begin{subfigure}{.3\textwidth}
		\centering
		\includegraphics[width=.8\linewidth]{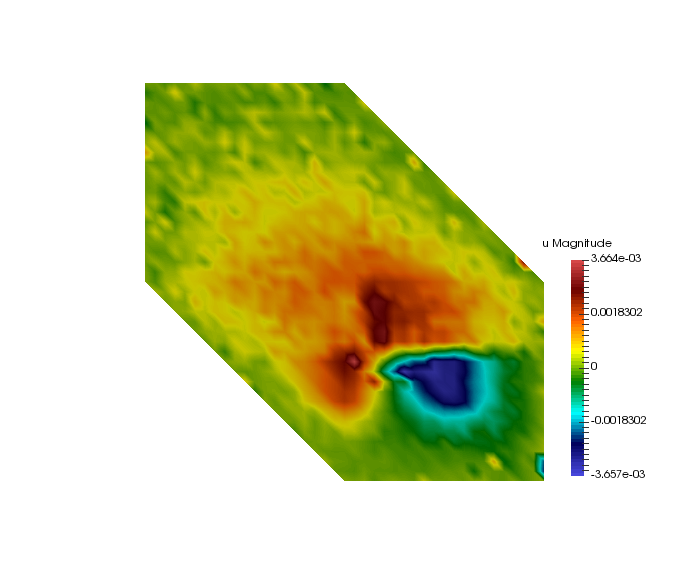}
		\caption{Difference $\phi_2^{\omega_3} - \phi_2^{\omega_3,\text{pred}}$}
	\end{subfigure}
	\caption{Exact multiscale basis functions $\phi_m^{\omega_3}$ and predicted multiscale basis functions $\phi_m^{\omega_3,\text{pred}}$ in the coarse neighborhood $\omega_3$ in Experiment 2.}
	\label{fig:98_basis_rand}
\end{figure}

Figures~\ref{fig:88_basis_rand}--\ref{fig:98_basis_rand} 
show the comparison of the multiscale basis functions 
in 2 respective coarse neighborhoods. 
It can be observed that the predicted multiscale basis functions 
are in good agreement with the exact ones. 
In particular, the neural network successfully interprets the 
high conductivity regions as the support localization feature 
of the multiscale basis functions. 
Tables~\ref{tab:exp2b1e}--\ref{tab:exp2me} record 
the mean error of the prediction by the neural networks, 
measured in the defined metric. 
Again, it can be seen that the prediction are of high accuracy. 
Table~\ref{tab:exp2se} records the mean error between the multiscale solution 
using the neural-network-based multiscale solver and using exact GMsFEM. 
we obtain a good approximation of the multiscale solution 
compared with the exact GMsFEM solver.

\begin{table}[ht!]
\centering
\begin{tabular}{c|cc|cc|cc}
Basis & \multicolumn{2}{|c|}{$\omega_1$} & \multicolumn{2}{|c|}{$\omega_2$} & \multicolumn{2}{|c}{$\omega_3$} \\
$m$ & $e_{L^2}$ & $e_{H^1}$ & $e_{L^2}$ & $e_{H^1}$ & $e_{L^2}$ & $e_{H^1}$ \\
\hline
1 & 0.55 & 0.91 & 0.37 & 3.02 & 0.20 & 0.63 \\
2 & 0.80 & 1.48 & 2.17 & 3.55 & 0.27 & 1.51 \\
\end{tabular}
\caption{Mean percentage error of multiscale basis functions $\phi_m^{\omega_i}$ in Experiment 2.}
\label{tab:exp2b1e}
\end{table}

\begin{table}[ht!]
\centering
\begin{tabular}{c|ccc}
& $e_{\ell^2}$   &  $e_{\ell^{\infty}}$  & $e_{F}$\\
\hline
Mean & 0.75 & 0.72 & 0.80
\end{tabular}
\caption{Percentage error of the local stiffness matrix $A_c^{K_0}$ in Experiment 2.}
\label{tab:exp2me}%
\end{table}

\begin{table}[ht!]
\centering
\begin{tabular}{c|cc}
& $e_{L^2}$  & $e_{a}$\\
\hline
Mean & 0.03 & 0.26
\end{tabular}
\caption{Percentage error of multiscale solution $u_{ms}$ in Experiment 2.}
\label{tab:exp2se}%
\end{table}

\section{Conclusion}
\label{sec:conclusion}
In this paper, we develop a method using deep learning techniques 
for fast computation of GMsFEM discretizations. 
Given a particular permeability field, the main ingredients of GMsFEM, 
including the multiscale basis functions and coarse-scale matrices, 
are computed in an offline stage by solving local problems. 
However, when one is interested in calculating GMsFEM discretizations 
for multiple choices of permeability fields, 
repeatedly formulating and solving such local problems could become 
computationally expensive or even infeasible. 
Multi-layer networks are used to represent the nonlinear mapping 
from the fine-scale permeability field coefficients 
to the multiscale basis functions and the coarse-scale parameters. 
The networks provide a direct fast approximation of the GMsFEM ingredients 
in a local neigorhood for any online permeability fields, 
in contrast to repeatedly formulating and solving local problems. 
Numerical results are presented to show the performance of our proposed method. 
We see that, given sufficient samples of GMsFEM discretizations for supervised training, 
deep neural networks are capable of providing 
reasonably close approximations of the exact GMsFEM discretization. 
Moreover, the small consistency error provides good approximations 
of multiscale solutions.

\bibliographystyle{plain}
\bibliography{references}

\end{document}